# Set-polynomials and polynomial extension of the Hales-Jewett Theorem

By V. BERGELSON and A. LEIBMAN*

### Abstract

An abstract, Hales-Jewett type extension of the polynomial van der Waerden Theorem [BL] is established:

THEOREM.   *Let $r, d, q \in \mathbb{N}$. There exists $N \in \mathbb{N}$ such that for any $r$-coloring of the set of subsets of $V = \{1, \ldots, N\}^d \times \{1, \ldots, q\}$ there exist a set $a \subset V$ and a nonempty set $\gamma \subseteq \{1, \ldots, N\}$ such that $a \cap (\gamma^d \times \{1, \ldots, q\}) = \emptyset$, and the subsets $a$, $a \cup (\gamma^d \times \{1\})$, $a \cup (\gamma^d \times \{2\})$, ..., $a \cup (\gamma^d \times \{q\})$ are all of the same color.*

This "polynomial" Hales-Jewett theorem contains refinements of many combinatorial facts as special cases. The proof is achieved by introducing and developing the apparatus of *set-polynomials* (polynomials whose coefficients are finite sets) and applying the methods of topological dynamics.

## 0. Introduction

0.1.  The celebrated van der Waerden Theorem ([W]) states that if the integers are partitioned into finitely many classes then one of the cells of the partition contains arbitrarily long arithmetic progressions.

In 1963, Hales and Jewett ([HJ]) published a far reaching generalization of van der Waerden's theorem. To formulate their result, we introduce some definitions. Let $A$ be a finite set, $A = \{s_1, \ldots, s_q\}$, and let $A^n$ denote the set of $n$-tuples with elements from $A$. A set $\{w^1, \ldots, w^q\} \subset A^n$, consisting of $q$ $n$-tuples $w^l = (w_1^l, \ldots, w_n^l)$, $l = 1, \ldots, q$, is *a combinatorial line*, if there exists a partition $\{1, \ldots, n\} = I \cup J$, $I \cap J = \emptyset$, $I \neq \emptyset$ such that $w_i^1 = \ldots = w_i^q$ for $i \in J$ and $w_i^l = s_l$ for $i \in I$ and $l = 1, \ldots, q$. Another way to describe a combinatorial line is the following. Let $t$ be a variable; consider the set of the words of length $n$ from the alphabet $A \cup \{t\}$. If $w(t)$ is a word in which the variable $t$ actually occurs, then the set $\{w(s_1), \ldots, w(s_q)\}$ (obtained by substituting the elements of $A$ instead of $t$) is a combinatorial line.

---
*The authors gratefully acknowledge support received from the National Science Foundation via grants DMS-9401093 and DMS-9706057.

0.2.  We are now in position to formulate the Hales-Jewett Theorem.

THEOREM HJ0 ([HJ]).  *For any $r \in \mathbb{N}$ there exists $N = N(q,r) \in \mathbb{N}$ such that if $A^N$ is partitioned into $r$ classes, then at least one of these classes contains a combinatorial line.*

If one takes $A = \{0, 1, \ldots, q-1\}$ and interprets the set $A^n$ as the set of base $q$ expansions of nonnegative integers smaller than $q^n$, then a combinatorial line forms an arithmetic progression (with difference $d = \sum_{i=0}^{n-1} a_i q^i$, where $a_i = 0$ or 1). This gives the van der Waerden Theorem, alluded to above.

On the other hand, if one takes $A$ to be a finite field, then $A^n$ may be viewed as a $n$-dimensional vector space over $A$. Then a combinatorial line is a line in the geometric sense in $A^n$ (i.e. an affine linear one-dimensional subspace of $A^n$).

0.3.  An interesting feature of HJ0 is that, as R. Graham ([G, p. 17]) has put it, "all the higher-dimensional analogs of Hales-Jewett [theorem] follow immediately from the 1-dimensional case." Let $t_1, \ldots, t_m$ be $m$ variables; consider the words over the alphabet $A \cup \{t_1, \ldots, t_m\}$. If $w(t_1, \ldots, t_m)$ is a word in which all the variables actually occur, then the set $\{w(s_{i_1}, \ldots, s_{i_m}) : s_{i_j} \in A, \ j = 1, \ldots, m\}$ will be called *a combinatorial $m$-space*. It follows from HJ0 that for any $r, m \in \mathbb{N}$ there exists $N = N(q, r, m)$ such that if $A^N$ is partitioned into $r$ classes then at least one of these classes contains a combinatorial $m$-space. To see this, note that a combinatorial line in $B^n$ for $B = A^m$ is a combinatorial $m$-space in $A^{mn}$.

Using the multidimensional or, rather, multiparameter version of HJ0 one derives from it the multidimensional van der Waerden Theorem as well as the Graham-Leeb-Rothschild geometric Ramsey Theorem ([GLR]).

0.4.  Given $r \in \mathbb{N}$, we define *an $r$-coloring* of a set $W$ as an $r$-fold partition of $W$ or, equivalently, as a mapping $\chi: W \longrightarrow \{1, \ldots, r\}$. We will say that a subset $\mathcal{M} \subseteq W$ *is monochromatic* (with respect to a coloring $\chi$), or elements of $\mathcal{M}$ are *of the same color* if $\chi$ is constant on $\mathcal{M}$.

Given an arbitrary set $W$, $\mathcal{F}(W)$ will denote the set of finite subsets of $W$ (including the empty one).

0.5.  We bring now two more formulations of the Hales-Jewett Theorem. These are the formulations which serve as the basis for our polynomial generalization of the theorem.

THEOREM HJ. *Let $q, r \in \mathbb{N}$. There exists $N \in \mathbb{N}$ such that for any $r$-coloring of $\mathcal{F}(\mathbb{N})^q$ (the set of $q$-tuples with entries from $\mathcal{F}(\mathbb{N})$) there exist a nonempty set $\gamma \subseteq \{1, \ldots, N\}$ and sets $a_1, \ldots, a_q \subseteq \{1, \ldots, N\}$ such that $\gamma$ is disjoint from $a_i$, $i = 1, \ldots, q$, and such that the set*

$$\Big\{(a_1, a_2, a_3, \ldots, a_q), (a_1 \cup \gamma, a_2, a_3, \ldots, a_q), (a_1, a_2 \cup \gamma, a_3, \ldots, a_q), \ldots,$$
$$(a_1, a_2, a_3, \ldots, a_q \cup \gamma)\Big\}$$

*is monochromatic.*

0.6. Given a set $W$, *an action* of $\mathcal{F}(W)$ on a topological space $X$ is a mapping $T$ from $\mathcal{F}(W)$ into the set of continuous self-mappings of $X$, $a \mapsto T^a$, satisfying the following condition: for any $a, b \in \mathcal{F}(W)$ with $a \cap b = \emptyset$ one has $T^{a \cup b} = T^a T^b$.

Notice that, for a set $W$ and $q \in \mathbb{N}$, $\mathcal{F}(W)^q$ is in a natural one-to-one correspondence with $\mathcal{F}(W \times \{1, \ldots, q\})$: the element $(a_1, \ldots, a_q)$ of $\mathcal{F}(W)^q$ corresponds to the element $(a_1 \times \{1\}) \cup \ldots \cup (a_q \times \{q\})$ of $\mathcal{F}(W \times \{1, \ldots, q\})$. We define an action of $\mathcal{F}(W)^q$ on $X$ as the action of $\mathcal{F}(W \times \{1, \ldots, q\})$.

THEOREM HJt. *Let $(X, \rho)$ be a compact metric space, let $q \in \mathbb{N}$ and let $T$ be an action of $\mathcal{F}(\mathbb{N})^q$ on $X$. Then for any $\varepsilon > 0$ there exists $N \in \mathbb{N}$ such that for any $x \in X$ there exist a nonempty set $\gamma \subseteq \{1, \ldots, N\}$ and sets $a_1, \ldots, a_q \subset \{1, \ldots, N\}$ such that $a_i \cap \gamma = \emptyset$, $i = 1, \ldots, q$, and*

$$\rho(T^{(a_1 \cup \gamma, a_2, \ldots, a_q)} x, T^{(a_1, a_2, \ldots, a_q)} x) < \varepsilon,$$
$$\rho(T^{(a_1, a_2 \cup \gamma, \ldots, a_q)} x, T^{(a_1, a_2, \ldots, a_q)} x) < \varepsilon,$$
$$\vdots$$
$$\rho(T^{(a_1, a_2, \ldots, a_q \cup \gamma)} x, T^{(a_1, a_2, \ldots, a_q)} x) < \varepsilon.$$

0.7 *Remark.* Being a general combinatorial fact about sets, the Hales-Jewett Theorem does not appeal to any special properties of $\mathbb{N}$, besides its cardinality. In the sections to follow Hales-Jewett type theorems will be usually formulated for a general infinite set $S$.

0.8. In the course of the proof of the Polynomial Szemerédi Theorem in [BL], a generalization of the van der Waerden Theorem in a different, polynomial direction was obtained:

THEOREM PvdW ([BL, Cor. 1.11]). *For any natural numbers $r$, $k$, $t$ and $m$, for any polynomials $p_{1,1}, \ldots, p_{1,t}, p_{2,1}, \ldots, p_{2,t}, \ldots, p_{k,1}, \ldots, p_{k,t}$ having zero constant terms and taking on integer values at the integers, for any vectors $v_1, \ldots, v_t \in \mathbb{Z}^m$ and for any $r$-coloring of $\mathbb{Z}^m$ there exist $a \in \mathbb{Z}^m$, $n \in \mathbb{Z}$, $n \neq 0$,*

*such that the set*

$$S(a,n) = \left\{a + \sum_{j=1}^{t} p_{i,j}(n)v_j, \ i = 1,\ldots,k\right\}$$

*is monochromatic.*

Or, in an equivalent, invariant form:

*For any natural numbers $r$, $t$ and $m$, any polynomial mapping $P\colon \mathbb{Z}^t \longrightarrow \mathbb{Z}^m$ with $P(0) = 0$, any finite set $F \subset \mathbb{Z}^t$ and any $r$-coloring of $\mathbb{Z}^m$ there are $a \in \mathbb{Z}^m$ and nonzero $n \in \mathbb{Z}$ such that the set $a + P(nF)$ is monochromatic.*

It follows, for example, that for any $k \in \mathbb{N}$ and any finite coloring of integers there always exist monochromatic configurations of the form $\{a, a+n^2, a+2n^2, \ldots, a+kn^2\}$ or, say, of the form $\{a, a+n, a+n^2, \ldots, a+n^k\}$.

Note that the original van der Waerden's Theorem corresponds to the case $t = m = 1$, $p_{i,1}(n) = in$, $i = 1,\ldots,k$.

0.9. *Remark.* In accordance with the general philosophy of Ramsey Theory (see [B2]), one should expect to get many polynomial configurations in one color. This is indeed so: one can show that the set of $n$ for which there is $a$ such that $S(a,n)$ is monochromatic (as well as the set of $n$ for which there is $a$ such that $a + P(nF)$ is monochromatic) is an IP*-set (see definition in 0.14 below; see also 7.14).

0.10. In this paper we obtain the polynomial Hales-Jewett Theorem (PHJ), which extends the Polynomial van der Waerden Theorem above in the way similar to that in which the Hales-Jewett Theorem generalizes the classical van der Waerden Theorem. Here is our main result, the Polynomial Hales-Jewett Theorem:

THEOREM PHJ. *For any $r,d,q \in \mathbb{N}$ there exists $N = N(r,d,q) \in \mathbb{N}$ such that for $V = \{1,\ldots,N\}^d \times \{1,\ldots,q\}$ the following holds: for any $r$-coloring of $\mathcal{F}(V)$ there exist $a \subset V$ and a nonempty set $\gamma \subseteq \{1,\ldots,N\}$ such that $a \cap (\gamma^d \times \{1,\ldots,q\}) = \emptyset$ and the sets*

$$a, \ a \cup (\gamma^d \times \{1\}), \ a \cup (\gamma^d \times \{2\}), \ \ldots, \ a \cup (\gamma^d \times \{q\})$$

*are all of the same color.*

0.11. The polynomial nature of Theorem PHJ is manifested by the integer $d$ in its formulation. This $d$ corresponds to "the degree" of polynomials that Theorem PHJ deals with. In particular, the "linear" case, $d = 1$, gives Theorem HJ. For a fixed $d$, the algebraic corollaries, which one obtains from Theorem PHJ, involve polynomials of degree not exceeding $d$. As an illustration, let us show how Theorem PHJ implies the following particular case of Polynomial van der Waerden Theorem PvdW:

> For any finite coloring of $\mathbb{Z}$ and any $d \in \mathbb{N}$ there are arbitrarily long monochromatic arithmetic progressions whose difference is of the form $n^d$ with $n > 0$.

Indeed, given $q \in \mathbb{N}$ and an $r$-coloring $\chi$ of $\mathbb{Z}$, introduce a coloring $\chi'$ on the set of subsets of $V = \{1,\ldots,N\}^d \times \{1,\ldots,q\}$ by $\chi'(a) = \chi(|a_1|+2|a_2|+\ldots+q|a_q|)$, where $|c|$ denotes the cardinality of set $c$ and $a_i = a \cap \bigl(\{1,\ldots,N\}^d \times \{i\}\bigr)$, $i = 1,\ldots,q$. Find $a$ and $\gamma$ as guaranteed by Theorem PHJ. Then the $(q+1)$-term arithmetic progression with base $|a_1| + 2|a_2| + \ldots + q|a_q|$ and difference $|\gamma^d| = |\gamma|^d$ is monochromatic.

The set-theoretical expressions of the form $\gamma^d \times \{l\}$ appearing in the formulation of PHJ are "polynomials" of quite special form, which correspond, roughly, to usual algebraic monomials of the form $lx^d$. To formulate a more general (at least in appearance) version of PHJ which would deal with more general set-theoretical expressions, one needs first to introduce *set-polynomials*. This is done in Section 2. The form of PHJ which uses the language of set-polynomials is given in Section 3 (Theorem 3.5). It can be shown, however, that the version of PHJ above is just as general.

0.12. The proof of Theorem PHJ uses ideas and methods from topological dynamics (which go back to the seminal paper of Furstenberg and Weiss [FW]). Here is an equivalent formulation of PHJ in the topological language (cf. Theorem HJt):

THEOREM PHJt. *Let $(X,\rho)$ be a compact metric space. For any $d$, $q \in \mathbb{N}$ there exists $N \in \mathbb{N}$ with the following property. Let $V = \{1,\ldots,N\}^d \times \{1,\ldots,q\}$ and let $T$ be an action of $\mathcal{F}(V)$ on $X$. Then for any $x \in X$ and any $\varepsilon > 0$ there exist $a \in \mathcal{F}(V)$ and a nonempty set $\gamma \subseteq \{1,\ldots,N\}$ such that $a \cap \bigl(\gamma^d \times \{1,\ldots,q\}\bigr) = \emptyset$ and $\rho(T^{a\cup(\gamma^d \times \{i\})}x, T^a x) < \varepsilon$ for every $i = 1,\ldots,q$.*

The equivalence of Theorems PHJ and PHJt suggests that in dealing with Ramsey-theoretical questions one has at his disposal not only the conventional language of combinatorics but also an equivalent language or, rather, method of topological dynamics. Whereas the combinatorial problems discussed in this paper are often formulated in a "finitary" way (i.e. as statements about finite sets), the equivalent topological-dynamical statements deal with potentially infinite compact spaces and continuous mappings thereof. It is the convenience of dealing with these "ideal" objects which is the reason for the effectiveness of methods of topological dynamics. For more information the reader is referred to [FW], [F], [FK], [BH], and [BBH].

A proof of the equivalence of "chromatic" and topological versions of the general form of PHJ is provided by Proposition 3.3 below. We want to add a few general remarks, explaining the nature of the equivalence of chromatic and topological statements similar to Theorem PHJ and Theorem PHJt, respectively. First of all, note that to derive the combinatorial proposition from a topological one, one considers the (compact) space of all $r$-colorings of the set in question and applies to it a general topological recurrence theorem. On the other hand, to derive a topological recurrence theorem from a combinatorial statement about monochromatic configurations which one always finds in any finite coloring, one uses the basic property of compact spaces, namely the existence of finite covers consisting of sets of arbitrarily small diameter, and, after inducing a finite coloring on the acting (semi) group with the help of such a cover, utilizes the combinatorial statement. We will allow ourselves to pass freely from the combinatorial language to the topological one and back when describing and deriving examples and corollaries of PHJ.

The proof of Theorem PHJ resembles that of the Polynomial van der Waerden Theorem in [BL]. However it has a certain novelty, which allows us to give a somewhat different proof even to the (linear) van der Waerden and Hales-Jewett Theorems (see subsection 1.1 below). The quintessence of this difference in the approach stems from the fact that while in [FW] (and [BL]) one reduces the setup to that of *a minimal* dynamical system, it is no longer possible to confine oneself to minimal systems when dealing with PHJ.

0.13.  We conclude this introductory section by collecting several examples of applications of PHJ. Numerous applications of PHJ include, for example, a strengthening and refinement of Theorem PvdW. To formulate it, we need the notion of IP-set. Let $G$ be a commutative semigroup. Given an infinite sequence $\{n_i\}_{i\in\mathbb{N}}$ in $G$ *the* IP-*set generated by* $\{n_i\}_{i\in\mathbb{N}}$ is the set of all finite sums of distinct elements from $\{n_i\}_{i\in\mathbb{N}}$. In other words, the IP-set generated by $\{n_i\}_{i\in\mathbb{N}}$ is the set $\{n_\alpha\}$ where $\alpha$ runs through the set of all nonempty finite subsets of $\mathbb{N}$ (that is, $\mathcal{F}(\mathbb{N})\setminus\{\emptyset\}$) and $n_\alpha := \sum_{i\in\alpha} n_i$. In particular, we have *additive* IP-*sets* in $(\mathbb{Z},+)$ and *multiplicative* IP-*sets* in the multiplicative positive integers $(\mathbb{N},\cdot)$.

THEOREM. *For any natural numbers $t$ and $m$, any polynomial mapping $P\colon\mathbb{Z}^t \longrightarrow \mathbb{Z}^m$ with $P(0) = 0$, any finite set $F \subset \mathbb{Z}^t$, any IP-set $\{n_\alpha = (n_{1,\alpha},\ldots,n_{t,\alpha}) : \alpha \subset \mathbb{N},\ 0 < |\alpha| < \infty\} \subset \mathbb{Z}^t$ and any finite coloring of $\mathbb{Z}^m$ there exist a vector $a \in \mathbb{Z}^m$ and a finite nonempty set $\gamma \subset \mathbb{N}$ such that the set*

$$\{a + P(n_{1,\gamma}v_1,\ldots,n_{t,\gamma}v_t) : (v_1,\ldots,v_t) \in F\}$$

*is monochromatic.*

(See Corollary 7.14 for a more general version of this proposition.)

0.14. A subset of a commutative semigroup $(G, +)$ is called *an* IP*-*set* if it has nontrivial intersection with any IP-subset of $G$. One can show that if $E \subseteq G$ is an IP*-set, then it is *syndetic*, namely there exists a finite set $F \subseteq G$ such that $\bigcup_{x \in F}(E - x) = G$ (where, by definition, $g \in E - x$ if and only if $g + x \in E$). In particular, when $G$ is a group, a set is syndetic if finitely many of its shifts cover $G$.

Since IP*-sets are syndetic, to say that a set $E \subseteq G$ is IP* means that it is, in a certain sense, quite large. The following remarks show that IP*-sets are even larger than they may look at the first sight. By Hindman's theorem ([H]), if $E \subseteq G$ is an IP*-set then for any IP-set $I \subseteq G$ the set $E \cap I$ contains an IP-set. This, in its turn, implies that the intersection of any finite family of IP*-sets is an IP*-set.

We can now reformulate Theorem 0.13 in the following way, which makes it evident that it extends and refines Theorem PvdW (see also Remark 0.9):

THEOREM.   *For any $t, m \in \mathbb{N}$, any polynomial mapping $P: \mathbb{Z}^t \longrightarrow \mathbb{Z}^m$ satisfying $P(0) = 0$, any finite set $F \subset \mathbb{Z}^t$ and any finite coloring $\chi$ of $\mathbb{Z}^m$, $\chi: \mathbb{Z}^m \longrightarrow \{1, \ldots, r\}$, there is $l$, $1 \le l \le r$, such that the set*

$$\Big\{(n_1, \ldots, n_t) : \text{there is } a \in \mathbb{Z}^m \text{ such that } \chi\big(a + P(n_1 v_1, \ldots, n_t v_t)\big) = l$$
$$\text{for all } (v_1, \ldots, v_t) \in F\Big\}$$

*is an* IP*-*subset of $\mathbb{Z}^t$.*

0.15. The following proposition corresponds to the special case of Theorem 0.14 where $t = 2q$, $m = 1$, $P(n_1, k_1, \ldots, n_q, k_q) = \sum_{i=1}^{q} n_i k_i$ and $F = \{(1, 1, 0, \ldots, 0), (0, 0, 1, 1, 0, \ldots, 0), \ldots, (0, 0, \ldots, 0, 1, 1)\}$.

PROPOSITION.   *Let $\{n_{1,i}\}_{i=1}^{\infty}, \{k_{1,i}\}_{i=1}^{\infty}, \ldots, \{n_{q,i}\}_{i=1}^{\infty}, \{k_{q,i}\}_{i=1}^{\infty}$ be sequences in $\mathbb{Z}$ and let $\{n_{1,\alpha}\}, \{k_{1,\alpha}\}, \ldots, \{n_{q,\alpha}\}, \{k_{q,\alpha}\}$ be the additive IP-sets generated by these sequences. Then for any finite coloring of $\mathbb{Z}$ there exists a monochromatic set of the form*

$$\{a, a + n_{1,\gamma} k_{1,\gamma}, \ldots, a + n_{q,\gamma} k_{q,\gamma}\}$$

*for some $a \in \mathbb{N}$ and a finite nonempty $\gamma \subset \mathbb{N}$.*

0.16. The possibility of diverse applications of PHJ stems from the fact that PHJ employs the basic set-theoretical ring-like operations of taking the (disjoint) unions and Cartesian products. Thus PHJ implies polynomial van der Waerden type theorems not only in $\mathbb{Z}$ but in general rings and commutative semigroups as well. Here is, for example, a multiplicative analogue of Proposition 0.15:

PROPOSITION. Let $\{\sigma_{1,i}\}_{i=1}^{\infty}, \ldots, \{\sigma_{q,i}\}_{i=1}^{\infty}, \{\pi_{1,i}\}_{i=1}^{\infty}, \ldots, \{\pi_{q,i}\}_{i=1}^{\infty}$ be sequences of natural numbers and let $\{\sigma_{1,\alpha}\}, \ldots, \{\sigma_{q,\alpha}\}$, and $\{\pi_{1,\alpha}\}, \ldots, \{\pi_{q,\alpha}\}$, for finite $\alpha \subset \mathbb{N}$, be, respectively, the additive and the multiplicative IP-sets generated by these sequences. Then for any finite coloring of $\mathbb{N}$ there exists a monochromatic set of the form

$$\{b, b\pi_{1,\gamma}^{\sigma_{1,\gamma}}, \ldots, b\pi_{q,\gamma}^{\sigma_{q,\gamma}}\}$$

for some $b \in \mathbb{N}$ and a finite nonempty $\gamma \subset \mathbb{N}$.

0.17. The following corollary of PHJ is an extension of Theorem PvdW to general fields (see 7.14 for a more general statement). Note that in the case of fields of finite characteristic the result provides a polynomial version of the geometric Ramsey theorem alluded to in 0.3.

THEOREM. Let $W$ and $V$ be vector spaces over an infinite field $\mathbf{K}$, let $P \colon W \longrightarrow V$ be a polynomial mapping with $P(0) = 0$ and let $F \subset W$ be a finite set. Then for any finite coloring $\chi \colon V \longrightarrow \{1, \ldots, r\}$ there are $a \in V$ and nonzero $n \in \mathbf{K}$ such that the set $a + P(nF)$ is monochromatic. Moreover, for some color $l$, $1 \leq l \leq r$, the set

$$\{n \in \mathbf{K} : \text{ there is } a \in V \text{ such that } \chi(a + P(nv)) = l \text{ for all } v \in F\}$$

is an IP*-set.

0.18. As a matter of fact, the class of mappings to which PHJ naturally applies is wider than the algebraic polynomials. Let $G$, $G'$ be commutative groups. We say that $f \colon G \longrightarrow G'$ is *a polynomial mapping* if $f$ trivializes after finitely many applications of the difference operator $(D_a f)(b) = f(a+b) - f(b)$ where $a \in G$. In other words, $f$ is a polynomial mapping if for some $k \in \mathbb{N}$ and any $a_1, \ldots, a_k \in G$ one has $D_{a_k} \ldots D_{a_1} f \equiv \mathbf{1}_{G'}$. In case of self-mappings of $\mathbb{Z}$ every polynomial mapping is a conventional algebraic polynomial (but not necessarily a member of $\mathbb{Z}[x]$: take, for example, $\frac{x^2-x}{2}$). On the other hand, for infinite fields of finite characteristic (or indeed, for any field which is an infinite dimensional extension of its prime subfield) one has an ample supply of polynomial mappings which are not algebraic polynomials. (For example, let $\mathbf{K}$ be an infinite field of characteristic 2. Let $E$ be a basis of $\mathbf{K}$, viewed as a vector space over $\mathbb{Z}_2$. The reader can check that for $k > 2$, the mapping $f \colon \mathbf{K} \longrightarrow \mathbf{K}$ defined by $f(a) = \binom{n(a)}{k}$, where $n(a)$ is the number of nonzero coordinates of $a \in \mathbf{K}$, is a polynomial mapping which is not a conventional polynomial.) It turns out that, for example, Theorem 0.17 holds true if the polynomial $P \colon W \longrightarrow V$ is replaced by a general polynomial mapping (see subsection 8.8 for a more general proposition).

0.19. The structure of the paper is as follows. Section 1 is devoted to a detailed presentation of two special cases of PHJ: the "linear" HJ Theorem and the simplest nontrivial case of PHJ. This section is included to let the reader get the spirit of our approach and stress the similarity between the number-theoretical polynomial van der Waerden theorem and its set-theoretical analogue, PHJ. However, an ambitious reader can safely skip it and go directly to Section 2.

In Sections 2 and 3 we give general definitions and establish some notation. With the help of this notation we formulate Theorem PHJ in another form (Theorem 3.4), which is more susceptible to the inductive procedure which is utilized for its proof. Sections 4 and 5 are devoted to some preliminary technical lemmas. Section 6 contains the proof of Theorem 3.4.

In Section 7 we develop a systematic procedure for derivation of corollaries from the PHJ Theorem. We use only the formulation of PHJ given in 0.10, so that the reader who is interested in applications but not in the proof of PHJ may see the strength and variety of results derivable from PHJ. An alternative possibility, which we decided not to pursue, was to derive the applications from an equivalent form of PHJ, which is brought in Section 3. This could make the derivation of applications somewhat more natural, but would require acquaintance with the material of Sections 2 and 3.

At the end of the paper we have placed an appendix, devoted to an (even more) abstract version of the PHJ Theorem and to discussion of the abstract polynomiality.

0.20 *Acknowledgment.* We thank H. Furstenberg, N. Hindman, Y. Katznelson and R. McCutcheon for helpful comments regarding an earlier draft of this paper.

## 1. Two special cases

To demonstrate the method of the proof of PHJ and to establish a certain parallelism between the proofs of the number-theoretical Polynomial van der Waerden Theorem and its abstract version, PHJ, we treat in this section two special cases of PHJ: the "linear" one, namely a variant of the classical Hales-Jewett Theorem, HJt, and the simplest nonlinear Hales-Jewett Theorem, which corresponds to $d = 2$, $q = 1$ of PHJt.

The reader may find instructive to compare the proofs of the somewhat modified version of (linear) the van der Waerden Theorem and that of Hales-Jewett which are given in two parallel columns having many identical portions. To ease the presentation and to emphasize the correspondence between the number-theoretical and set-theoretical notions, we will abide in this section by the following notational agreement: "+" will be used both for addition in $\mathbb{N}$ and for operation of taking (disjoint) unions of sets, "−" will be used not only for subtraction in $\mathbb{N}$ (when the minuend is not smaller than the subtrahend) but also instead of the set-theoretical difference "\" in expressions of the form $A \setminus B$ where $B \subseteq A$. The sign "·" will be used for both the multiplication in $\mathbb{N}$ and for the operation of taking the Cartesian products (and will often be omitted). The sign "$\preceq$" will mean either the usual inequality "$\leq$" or the set-theoretical containment "$\subseteq$". "0" will mean either zero or the empty set $\emptyset$.

1.1.   Let $(X, \rho)$ be a compact metric space. Let $q \in \mathbb{N}$.

Denote the set of nonnegative integers by $\mathbb{F}$. Let $T$ be a continuous self-mapping of $X$. Let $A$ be a set consisting of $q$ pairwise distinct natural numbers

$$A = \{p_i \in \mathbb{N} : i = 1, \ldots, q\};$$

assume without loss of generality that $p_1 < p_2 < \ldots < p_q$.

Let $S$ be an infinite set, denote $\mathcal{F}(S)$ by $\mathbb{F}$. Let $V = \{1, \ldots, q\} \times S$ and let $T$ be an action of $\mathcal{F}(V)$ on $X$. Put $p_i = \{1, \ldots, i\}$, $i = 1, \ldots, q$, and $A = \{p_1, \ldots, p_q\}$.

We are going to prove the following (two) proposition(s):

PROPOSITION L (vdW, HJ).   *For any $\varepsilon > 0$ there exists $N \in \mathbb{F}$, such that for any $x \in X$ there exist $n \preceq N$, $n \neq 0$, and $a \preceq p_q(N - n)$ such that for any $p \in A$,*

$$\rho(T^{a+pn}x, T^a x) < \varepsilon.$$

*Proof.* We will prove this proposition by induction on $q$. Define $B$ by

$$B = \{p_i - p_1, \ i = 2, \ldots, q\}.$$

Since $B$ contains $q - 1$ elements, we may assume that the statement to prove is valid for $B$, that is, for any $\varepsilon > 0$, there exists $N \in \mathbb{F}$ such that for any $x \in X$ there exist $n \preceq N$, $a \preceq (p_q - p_1)(N - n)$, such that $n \neq 0$ and for every $r \in B$ one has $\rho(T^{a+rn}x, T^a x) < \varepsilon$.

Let $\varepsilon > 0$. Let $k \in \mathbb{N}$ be such that among any $k + 1$ points of $X$ there are two points at a distance less than $\varepsilon/2$.

1.1.1. Put $\varepsilon_0 = \varepsilon/2k$. By the induction hypothesis, there exists $N_0 \in \mathbb{F}$ such that for any $x \in X$ there exist $n \preceq N_0$ and $a \preceq (p_q - p_1)(N_0 - n)$ such that $n \neq 0$ and for every $r \in B$ one has $\rho(T^{a+rn}x, T^a x) < \varepsilon_0$.

Let $\varepsilon_1 > 0$ be such that the inequality $\rho(y_1, y_2) < \varepsilon_1$ implies the inequality

$$\rho(T^b y_1, T^b y_2) < \varepsilon/2k$$

for any $b \preceq p_q N_0$. Let $N_1 \in \mathcal{F}$ be such that $N_1 \cap N_0 = 0$ (this disjointness condition concerns Proposition L(HJ) only) and for any $x \in X$ there exist $n \preceq N_1$ and $a \preceq (p_q - p_1)(N_1 - n)$ such that $n \neq 0$ and for every $r \in B$ one has $\rho(T^{a+rn}x, T^a x) < \varepsilon_1$.

Continue this process: assume that $\varepsilon_0, \ldots, \varepsilon_i$ and $N_0, \ldots, N_i \in \mathcal{F}$ have been already chosen. Let $\varepsilon_{i+1} > 0$ be such that the inequality $\rho(y_1, y_2) < \varepsilon_{i+1}$ implies the inequality

$$\rho(T^b y_1, T^b y_2) < \varepsilon/2k$$

for any $b \preceq p_q(N_0 + \ldots + N_i)$. Let $N_{i+1} \in \mathcal{F}$ be such that $N_{i+1} \cap (N_0 \cup \ldots \cup N_i) = 0$ (again, this disjointness condition is relevant for Proposition L(HJ) only) and for any $x \in X$ there exist $n \preceq N_{i+1}$ and $a \preceq (p_q - p_1)(N_{i+1} - n)$, such that $n \neq 0$ and for every $r \in B$ one has $\rho(T^{a+rn}x, T^a x) < \varepsilon_{i+1}$.

Continue the process of choosing $\varepsilon_i$, $N_i$ up to $i = k$, and put $N = N_0 + \ldots + N_k$.

1.1.2. Now fix an arbitrary point $x \in X$.

Applying the definition of $N_k$ to the point

$$y_k = T^{p_1 N_k} x,$$

we find $n_k \preceq N_k$, $n_k \neq 0$, and $a_k \preceq (p_q - p_1)(N_k - n_k)$ such that for every $r \in B$ we have

$$\rho(T^{rn_k + a_k} y_k, T^{a_k} y_k) < \varepsilon_k.$$

Then, applying the definition of $N_{k-1}$ to the point

$$y_{k-1} = T^{p_1(N_k + N_{k-1}) + a_k} x,$$

we find $n_{k-1} \preceq N_{k-1}$, $n_{k-1} \neq 0$, and $a_{k-1} \preceq (p_q - p_1)(N_{k-1} - n_{k-1})$ such that for every $r \in B$ we have

$$\rho(T^{rn_{k-1} + a_{k-1}} y_{k-1}, T^{a_{k-1}} y_{k-1}) < \varepsilon_{k-1}.$$

Continue this process: suppose that we have already found $n_k, \ldots, n_i$, $a_k, \ldots, a_i$. Applying the definition of $N_{i-1}$ to the point

$$y_{i-1} = T^{p_1(N_k + \ldots + N_{i-1}) + a_k + \ldots + a_i} x,$$

find $n_{i-1} \preceq N_{i-1}$, $n_{i-1} \neq 0$, and $a_{i-1} \preceq (p_q - p_1)(N_{i-1} - n_{i-1})$ such that for every $r \in B$ we have

$$\rho(T^{rn_{i-1} + a_{i-1}} y_{i-1}, T^{a_{i-1}} y_{i-1}) < \varepsilon_{i-1}.$$

Continue the process of choosing $n_i$, $a_i$ up to $i = 0$.

1.1.3. For every $0 \leq i \leq k$ we have $0 \neq n_i \preceq N_i$ and $a_i \preceq (p_q - p_1)(N_i - n_i)$;

therefore, for any $0 \leq i \leq j \leq k$ we have

(1.1)
$$p(n_j + \ldots + n_i) + a_j + \ldots + a_0$$
$$+ p_1(N_j + \ldots + N_0 - n_j - \ldots - n_0)$$
$$\preceq p_q(n_j + \ldots + n_{i+1})$$
$$+ (p_q - p_1)(N_j - n_j + \ldots + N_0 - n_0)$$
$$+ p_1(N_j + \ldots + N_0 - n_j - \ldots - n_0)$$
$$= p_q(N_0 + \ldots + N_j).$$

besides, $N_i \cap N_l = 0$ for $i \neq l$. Therefore, for any $0 \leq i \leq j \leq k$ we have

(1.2)
$$\Big( p(n_j + \ldots + n_i) + a_j + \ldots + a_0$$
$$+ p_1(N_j + \ldots + N_0 - n_j - \ldots - n_0) \Big)$$
$$\cap \Big( a_{j+1} + (p - p_1) n_{j+1} \Big) = 0.$$

And, for any $0 \leq j \leq k$,

(1.3)
$$a_k + \ldots + a_0 + p_1(N_k + \ldots + N_0 - n_j - \ldots - n_0)$$
$$\preceq (p_q - p_1)(N_k - n_k + \ldots + N_0 - n_0) + p_1(N_k + \ldots + N_0 - n_j - \ldots - n_0)$$
$$\preceq p_q(N_k + \ldots + N_0 - n_j - \ldots - n_0).$$

1.1.4. Define points $x_i$, $i = 0, \ldots, k$, by
$$x_i = T^{a_k + \ldots + a_0 + p_1(N_k + \ldots + N_0 - n_i - \ldots - n_0)} x.$$

We are going to show that for any $0 \leq i \leq j \leq k$ and any $p \in A$,

(1.4) $$\rho(T^{p(n_j + \ldots + n_{i+1})} x_j, x_i) \leq \frac{\varepsilon}{2k}(j - i).$$

We will prove this by induction on $j - i$; when $j = i$ the statement is trivial. We will derive the validity of (1.4) for $i, j$, where $i < j$, from its validity for $i, j - 1$.

*An important remark.* Note that the indices of the entries of expressions like "$n_j + \ldots + n_{i+1}$ for $i \leq j$" are assumed to be decreasing; when $i = j$ we assume that this expression is vacuous and equals zero. Such situations will often occur in our considerations, and we will always abide by this agreement.

By the definition of $n_j$,
$$\rho(T^{a_j + (p - p_1)n_j} y_j, T^{a_j} y_j) < \varepsilon_j,$$
where
$$y_j = T^{p_1(N_k + \ldots + N_j) + a_k + \ldots + a_{j+1}} x.$$

So, by the choice of $\varepsilon_j$ (and (1.1)),
$$\rho(T^{p(n_{j-1} + \ldots + n_{i+1}) + a_{j-1} + \ldots + a_0 + p_1(N_{j-1} + \ldots + N_0 - n_{j-1} - \ldots - n_0)} T^{a_j + (p - p_1)n_j} y_j,$$
$$T^{p(n_{j-1} + \ldots + n_{i+1}) + a_{j-1} + \ldots + a_0 + p_1(N_{j-1} + \ldots + N_0 - n_{j-1} - \ldots - n_0)} T^{a_j} y_j) < \varepsilon/2k.$$

Using the definition of $y_j$, $x_j$ (and (1.2)), we see

$$T^{p(n_{j-1}+\ldots+n_{i+1})+a_{j-1}+\ldots+a_0+p_1(N_{j-1}+\ldots+N_0-n_{j-1}-\ldots-n_0)}T^{a_j+(p-p_1)n_j}y_j$$
$$= T^{p(n_j+\ldots+n_{i+1})+a_k+\ldots+a_0+p_1(N_k+\ldots+N_0-n_j-\ldots-n_0)}x$$
$$= T^{p(n_j+\ldots+n_{i+1})}x_j$$

and

$$T^{p(n_{j-1}+\ldots+n_{i+1})+a_{j-1}+\ldots+a_0+p_1(N_{j-1}+\ldots+N_0-n_{j-1}-\ldots-n_0)}T^{a_j}y_j$$
$$= T^{p(n_{j-1}+\ldots+n_{i+1})+a_k+\ldots+a_0+p_1(N_k+\ldots+N_0-n_{j-1}-\ldots-n_0)}x$$
$$= T^{p(n_{j-1}+\ldots+n_{i+1})}x_{j-1}.$$

Since, by the induction hypothesis,

$$\rho(T^{p(n_{j-1}+\ldots+n_{i+1})}x_{j-1}, x_i) \leq \frac{\varepsilon}{2k}(j-i-1),$$

we obtain (1.4).

1.1.5. By the choice of $k$, among the $k+1$ points $x_0, \ldots, x_k$ there are two, say $x_i$, $x_j$, $0 \leq i < j \leq k$, for which $\rho(x_i, x_j) < \varepsilon/2$. Put

$$n = n_j + \ldots + n_{i+1},$$
$$a = a_k + \ldots + a_0 + p_1(N_k + \ldots + N_0 - n_j - \ldots - n_0).$$

Then $x_j = T^a x$ and

$$\rho(T^{a+pn}x, T^a x) = \rho(T^{pn}x_j, x_j)$$
$$\leq \rho(T^{pn}x_j, x_i) + \rho(x_j, x_i) < \varepsilon(j-i)/2k + \varepsilon/2 \leq \varepsilon.$$

Furthermore, $n \preceq N$, $n \neq 0$ and $a \preceq p_q(N-n)$ by (1.3). This proves Proposition L. □

1.2. We pass to the simplest nonlinear case of PHJ, the case of "a single square." When convenient, we will keep using the notational agreement made at the beginning of Section 1.

Let $(X, \rho)$ be a metric compact space, let $S$ be an infinite set and let $T$ be an action of $\mathcal{F}(S \times S)$ on $X$. We are going to prove the following proposition:

PROPOSITION Q. *For any $\varepsilon > 0$ there exists $N \in \mathcal{F}(S \times S)$ such that for any $x \in X$ there exist a nonempty $n \subseteq N$ and a set $a \subset N \times N$ such that $a \cap (n \times n) = \emptyset$ and*
$$\rho(T^{a+n\times n}x, T^a x) < \varepsilon.$$

*Remark.* The chromatic version of Proposition Q says that if $\mathcal{F}(S \times S)$ is finitely colored, $\mathcal{F}(S \times S) = \bigcup_{l=1}^{r} C_l$, then one of $C_l$ contains a pair $a, a \cup (n \times n)$ with $a \cap (n \times n) = \emptyset$. One of the consequences of this statement is the following (known) fact: for any finite coloring of $\mathbb{N}$ there exist $n \in \mathbb{N}$ and monochromatic $x$, $y$ such that $x - y = n^2$ (to see this induce the coloring of $\mathcal{F}(S \times S)$ by $\tilde{\chi}(a) = \chi(|a|)$).

1.3. The proof of Proposition Q uses still another form of linear HJ:

PROPOSITION L'. *Let $(X, \rho)$ be a compact metric space, let $S$ be an infinite set and let $T$ be an action of $\mathcal{F}(S \times S)$ on $X$. Let $A$ be a finite subset of $\mathcal{F}(S)$ and let $H$ be the union of its elements. Then, for any $\varepsilon > 0$, there exists $N \in \mathcal{F}(S)$ such that for any $x \in X$ there exist a nonempty set $n \subseteq N$ and a subset $a \subset H \times N \cup N \times H$ such that $a \cap (H \times n \cup n \times H) = \emptyset$ and for every $p \in A$ one has*

$$\rho(T^{a+(p \times n + n \times p)}x, T^a x) < \varepsilon.$$

Instead of deriving Proposition L' directly from Proposition L(HJ) ($\equiv$ Theorem HJt), we prefer to first reformulate it combinatorially and to derive this combinatorial version from Theorem HJ. A similar method will be utilized in the proof of the general PHJ.

PROPOSITION L''. *Let $S$ be an infinite set, let $A$ be a finite subset of $\mathcal{F}(S)$ and let $H$ be the union of the elements of $A$. Then, for any $r \in \mathbb{N}$ and any mapping $\chi \colon \mathcal{F}(S \times S) \longrightarrow \{1, \ldots, r\}$, there exists $N \in \mathcal{F}(S)$ such that there exist a nonempty set $n \subseteq N$ and a subset $a \subset H \times N \cup N \times H$ such that $a \cap (H \times n \cup n \times H) = \emptyset$ and $\chi$ is constant on the set*

$$\{p \times n \cup n \times p : p \in A\}.$$

*Proof.* Let $A = \{p_1, \ldots, p_q\}$. Let $W = \{1, \ldots, q\} \times \mathbb{N}$ and let $\sigma$ be an arbitrary embedding of $\mathbb{N}$ into $S \setminus H$. Define a mapping $\varphi \colon \mathcal{F}(W) \longrightarrow \mathcal{F}(S \times S)$ by the rule

$$\varphi(b) = \bigcup_{k \in \mathbb{N}} (p_{d_k} \times \{\sigma(k)\}) \cup (\{\sigma(k)\} \times p_{d_k}), \ b \in \mathcal{F}(W),$$

where $d_k = \#(b \cap (\{1, \ldots, q\} \times \{k\}))$. Now $\chi$ defines the finite coloring $\chi \circ \varphi$ of $\mathcal{F}(W)$, and the $\varphi$-images of the sets $(a_1, \ldots, a_q) \in \mathcal{F}(W)$ and $\gamma \subset \mathbb{N}$ given by Theorem HJ satisfy the conclusion of Proposition L''. □

*Proof of Proposition* Q. For $n, m \in \mathcal{F}(S)$ we will denote by $\mathrm{Sym}(n, m)$ the set $(n \times m) \cup (m \times n) \in \mathcal{F}(S \times S)$, and by $n^2$ the set $n \times n \in \mathcal{F}(S \times S)$. Note that $(n+m)^2 - n^2 - m^2 = \mathrm{Sym}(n, m)$.

Let $\varepsilon > 0$. Let $k \in \mathbb{N}$ be such that there are two points at a distance less than $\varepsilon/2$ among any $k+1$ points of $X$.

1.4.1. Let $N_0$ be an arbitrary finite nonempty subset of $S$, put $H_0 = N_0$, $A_0 = \emptyset$.

Let $\varepsilon_1 > 0$ be such that the inequality $\rho(y_1, y_2) < \varepsilon_1$ implies the inequality $\rho(T^b y_1, T^b y_2) < \varepsilon/2k$ for any $b \subseteq H_0^2$. Consider the system of sets

$$A_1 = \{p : p \subseteq H_0\}.$$

By Proposition L', there exists $N_1 \in \mathcal{F}(S)$, $N_1 \cap H_0 = \emptyset$, such that for any $x \in X$ there exist a nonempty set $n \subseteq N_1$ and a set $a \subset \operatorname{Sym}(H_0, N_1)$ such that $a \cap \operatorname{Sym}(H_0, n) = \emptyset$ and for every $p \in A_1$ one has

$$\rho(T^{a+\operatorname{Sym}(p,n)} x, T^a x) < \varepsilon_1.$$

Continue this process: assume that positive numbers $\varepsilon_0, \ldots, \varepsilon_i$ and finite sets $N_0, \ldots, N_i$ have been already chosen, and assume that the systems $A_0, \ldots, A_i$ have been already defined. Denote $H_i = N_0 + \ldots + N_i$. Let $\varepsilon_{i+1} > 0$ be such that the inequality $\rho(y_1, y_2) < \varepsilon_{i+1}$ implies the inequality $\rho(T^b y_1, T^b y_2) < \varepsilon/2k$ for all $b \subseteq H_i^2$. Consider the system

$$A_{i+1} = \{p \subseteq H_i\}.$$

By Proposition L', there exists $N_{i+1} \in \mathcal{F}(S)$, $N_{i+1} \cap H_i = \emptyset$, such that for any $x \in X$ there exist a nonempty set $n \subseteq N_{i+1}$ and a set $a \subset \operatorname{Sym}(H_i, N_{i+1})$ such that $a \cap \operatorname{Sym}(H_i, n) = \emptyset$ and for every $p \in A_{i+1}$ one has

$$\rho(T^{a+\operatorname{Sym}(p,n)} x, T^a x) < \varepsilon_{i+1}.$$

Continue the process of choosing $\varepsilon_i$, $N_i$ and defining $A_i$ up to $i = k$, and put $N = N_0 + \ldots + N_k$.

1.4.2. Fix now an arbitrary point $x \in X$.
Find, applying the definition of $N_k$ to the point

$$y_k = T^{N_k^2} x,$$

a nonempty set $n_k \subseteq N_k$ and a set $a_k \subseteq \operatorname{Sym}(H_{k-1}, N_k)$ such that $a_k \cap \operatorname{Sym}(H_{k-1}, n_k) = \emptyset$ and such that for every $p \in A_k$ one has

$$\rho(T^{\operatorname{Sym}(p,n_k)+a_k} y_k, T^{a_k} y_k) < \varepsilon_k.$$

Then, applying the definition of $N_{k-1}$ to the point

$$y_{k-1} = T^{N_k^2 + N_{k-1}^2 + a_k} x,$$

find a nonempty set $n_{k-1} \subseteq N_{k-1}$, and a set $a_{k-1} \subset \operatorname{Sym}(H_{k-2}, N_{k-1})$ such that $a_{k-1} \cap \operatorname{Sym}(H_{k-2}, n_{k-1}) = \emptyset$ and such that for every $p \in A_{k-1}$ one has

$$\rho(T^{\operatorname{Sym}(p,n_{k-1})+a_{k-1}} y_{k-1}, T^{a_{k-1}} y_{k-1}) < \varepsilon_{k-1}.$$

Continue this process: suppose that we have already found $n_k, \ldots, n_i$, $a_k, \ldots, a_i$. Applying the definition of $N_{i-1}$ to the point

$$y_{i-1} = T^{N_k^2 + \ldots + N_{i-1}^2 + a_k + \ldots + a_i} x,$$

find nonempty sets $n_{i-1} \subseteq N_{i-1}$, and sets $a_{i-1} \subset \operatorname{Sym}(H_{i-2}, N_{i-1})$ such that $a_{i-1} \cap \operatorname{Sym}(H_{i-2}, N_{i-1}) = \emptyset$ and such that for every $p \in A_{i-1}$ one has

$$\rho(T^{\operatorname{Sym}(p, n_{i-1}) + a_{i-1}} y_{i-1}, T^{a_{i-1}} y_{i-1}) < \varepsilon_{i-1}.$$

Continue the process of choosing $n_i$, $a_i$ up to $i = 0$.

1.4.3. The chosen sets $N_i$, $n_i$, $a_i$, $i = 0, \ldots, k$, satisfy the following conditions:

$$N_i \cap N_j = \emptyset, \ i \neq j, \ i, j = 0, \ldots, k,$$

$$n_i \subseteq N_i, \ i = 0, \ldots, k,$$

$$a_i \subset \operatorname{Sym}(N_0 + \ldots + N_{i-1}, N_i), \ i = 1, \ldots, k,$$

$$a_i \cap \operatorname{Sym}(N_0 + \ldots + N_{i-1}, n_i) = \emptyset, \ i = 1, \ldots, k.$$

It is clear, therefore, that

$$a_i \cap N_j^2 = \emptyset, \ i, j = 0, \ldots, k,$$

$$a_i \cap \operatorname{Sym}(n_j, n_l) = \emptyset, \ i, j, l = 0, \ldots, k,$$

$$a_i \cap a_j = \emptyset, \ i \neq j, \ i, j = 0, \ldots, k.$$

It is easy to establish from this that, for $0 \leq i < j \leq k$,

(1.5)
$$\bigl((n_{j-1} + \ldots + n_i)^2 + a_{j-1} + \ldots + a_0 + N_{j-1}^2 + \ldots + N_0^2 - n_{j-1}^2 - \ldots - n_0^2\bigr)$$
$$\cap \bigl(a_j + (n_j + \ldots + n_{i+1})^2 - (n_{j-1} + \ldots + n_{i+1})^2 - n_j^2\bigr) = \emptyset$$

and

(1.6) $\quad (n_j + \ldots + n_{i+1})^2 \cap \bigl(a_k + \ldots + a_0 + N_k^2 + \ldots + N_0^2 - n_j^2 - n_0^2\bigr) = \emptyset.$

1.4.4. Define points $x_i$, $i = 0, \ldots, k$, by

$$x_i = T^{a_k + \ldots + a_0 + N_k^2 + \ldots + N_0^2 - n_i^2 - \ldots - n_0^2} x.$$

We are going to show that, for any $0 \leq i \leq j \leq k$,

(1.7) $$\rho(T^{(n_j + \ldots + n_{i+1})^2} x_j, x_i) \leq \frac{\varepsilon}{2k}(j - i).$$

We will prove this by induction on $j - i$; when $j = i$ the statement is trivial. We will derive the validity of (1.7) for $i, j$, where $i < j$, from its validity for $i, j - 1$.

By the definition of $n_j$,
$$\rho(T^{a_j+(n_j+...+n_{i+1})^2-(n_{j-1}+...+n_{i+1})^2-n_j^2}y_j, T^{a_j}y_j) < \varepsilon_j,$$
where
$$y_j = T^{N_k^2+...+N_j^2+a_k+...+a_{j+1}}x.$$

So, by the choice of $\varepsilon_j$,
$$\rho(T^{(n_{j-1}+...+n_{i+1})^2+a_{j-1}+...+a_0+N_{j-1}^2+...+N_0^2-n_{j-1}^2-...-n_0^2}$$
$$T^{a_j+(n_j+...+n_{i+1})^2-(n_{j-1}+...+n_{i+1})^2-n_j^2}y_j,$$
$$T^{(n_{j-1}+...+n_{i+1})^2+a_{j-1}+...+a_0+N_{j-1}^2+...+N_0^2-n_{j-1}^2-...-n_0^2}T^{a_j}y_j) < \frac{\varepsilon}{2k}.$$

Using (1.5), (1.6) and the definition of $y_j$, $x_j$, we see
$$T^{(n_{j-1}+...+n_{i+1})^2+a_{j-1}+...+a_0+N_{j-1}^2+...+N_0^2-n_{j-1}^2-...-n_0^2}$$
$$T^{a_j+(n_j+...+n_{i+1})^2-(n_{j-1}+...+n_{i+1})^2-n_j^2}y_j$$
$$= T^{(n_j+...+n_{i+1})^2+a_k+...+a_0+N_k^2+...+N_0^2-n_j^2-...-n_0^2}x$$
$$= T^{(n_j+...+n_{i+1})^2}x_j$$

and
$$T^{(n_{j-1}+...+n_{i+1})^2+a_{j-1}+...+a_0+N_{j-1}^2+...+N_0^2-n_{j-1}^2-...-n_0^2}T^{a_j}y_j$$
$$= T^{(n_{j-1}+...+n_{i+1})^2+a_k+...+a_0+N_k^2+...+N_0^2-n_{j-1}^2-...-n_0^2}x$$
$$= T^{(n_{j-1}+...+n_{i+1})^2}x_{j-1}.$$

Since, by the induction hypothesis,
$$\rho(T^{(n_{j-1}+...+n_{i+1})^2}x_{j-1}, x_i) \leq \frac{\varepsilon}{2k}(j-i-1),$$
we obtain (1.7).

1.4.5. By the choice of $k$, among the $k+1$ points $x_0, \ldots, x_k$ there are two, say $x_i, x_j$, $0 \leq i < j \leq k$, for which $\rho(x_i, x_j) < \varepsilon/2$. Put
$$n = n_j + \ldots + n_{i+1}, \ a = a_k + \ldots + a_0 + N_k^2 + \ldots + N_0^2 - n_j^2 - \ldots - n_0^2.$$
Then $x_j = T^a x$; hence,
$$\rho(T^{a+n^2}x, T^a x) = \rho(T^{n^2}x_j, x_j) \leq \rho(T^{n^2}x_j, x_i) + \rho(x_j, x_i) < \frac{\varepsilon}{2k}(j-i) + \frac{\varepsilon}{2} \leq \varepsilon.$$

Furthermore, $n \subseteq N$, $n \neq \emptyset$ and $a \subseteq N^2 - n^2$ by (1.6). $\square$

## 2. Set-polynomials

In this section we define *set-polynomials* and derive some elementary facts about them. This formalism will be utilized in the following sections, where we state and prove another, more abstract form of Theorem PHJ.

2.1. Let us remind that for an arbitrary set $W$ we denote by $\mathcal{F}(W)$ the set of all finite subsets of $W$. Let $a, b \in \mathcal{F}(W)$. It will be convenient to denote the union $a \cup b$ by $a + b$. When $b \subseteq a$, the set-theoretical difference $a \setminus b$ will be written as $a - b$. The product, $ab$, of sets $a \in \mathcal{F}(W_1)$ and $b \in \mathcal{F}(W_2)$ is their Cartesian product $a \times b \in \mathcal{F}(W_1 \times W_2)$. Given a set $m$ and $d \in \mathbb{N}$, $m^d$ will stand for the $d$-fold Cartesian power of $m$. Note that if $m \in \mathcal{F}(W)$ then $m^d \in \mathcal{F}(W^d)$.

We will freely use both the "old" notation "$a \cup b$", "$a \setminus b$", "$a \times b$" and the "new" one "$a + b$", "$a - b$", "$ab$". The cardinality of $a$ will be denoted by either $\#a$ or $|a|$.

2.2. The stage is now set for the main definition of this section, namely, the definition of *the set-polynomial*. Just as the ordinary polynomials, say the polynomials of one variable over a ring, are defined via the operations of addition and multiplication, the set-polynomials are defined with the help of the operations of the addition (unions) and the (Cartesian) multiplication of finite sets.

For instance, a set-theoretical analogue of the polynomial $P(x) = x^2$ is the set-polynomial $P(n) = n^2$, where $n \in \mathcal{F}(W)$ and $n^2 = n \times n \in \mathcal{F}(W^2)$. A more general example of the "quadratic" set-polynomial is given by

$$P(n) = n^2 + b_1 n + n b_2 + c, \ n \in \mathcal{F}(W),$$

where $b_1, b_2 \in \mathcal{F}(W)$, $c \in \mathcal{F}(W^2)$. This last example does not exhaust all possible quadratic set-polynomials, because our definition depends on the dimension $d$ of the set $W^d$ from which coefficients are taken. The following set-polynomial is a set-polynomial of the second degree in $\mathcal{F}(W^3)$:

(2.1)   $a_1 \times n \times n + n \times a_2 \times n + n \times n \times a_3 + b_1 \times n + B_2(n) + n \times b_3 + c,$

where $a_1, a_2, a_3 \in \mathcal{F}(W)$, $b_1, b_2, b_3 \in \mathcal{F}(W^2)$, $B_2(n) = \{(s_1, s_2, s_3) : s_2 \in n, s_1 \times s_3 \in b_2\}$, and $c \in \mathcal{F}(W^3)$. More examples will be given after the formal definition, to which we now pass.

2.3. From now on $S$ will be an infinite set. Fix a positive integer $D$; the main object of our consideration will be $\mathcal{F}(V)$, where $V = S^D$.

Given a set of *coefficients* $\{P_\alpha\}_{\alpha \subseteq \{1,\ldots,D\}}$, $P_\alpha \in \mathcal{F}(S^{|\bar{\alpha}|})$ (where $\bar{\alpha} = \{1,\ldots,D\} \setminus \alpha$), the *set-polynomial* $P$ with coefficients $P_\alpha$ is the mapping

$P: \mathcal{F}(S) \longrightarrow \mathcal{F}(V)$ defined by
(2.2)
$$P(n) = \bigcup_{\alpha \subseteq \{1,\ldots,D\}} \{s = (s_1,\ldots,s_D) \in V : \ s_{\bar{\alpha}} \in P_\alpha \text{ and } s_i \in n \text{ for } i \in \alpha \ \},$$

where we denote $s_\beta = (s_{i_1},\ldots,s_{i_d})$ when $\beta = \{i_1,\ldots,i_d\}$ with $i_1 < \ldots < i_d$. The term of the union (2.2) which corresponds to $\alpha \subseteq \{1,\ldots,D\}$ will be called *the $\alpha$-term* of $P$. (A special agreement is needed for the case $\alpha = \{1,\ldots,D\}$, $\bar{\alpha} = \emptyset$. We assume $S^0 = \{\emptyset\}$, $\mathcal{F}(S^0) = \{\emptyset, \{\emptyset\}\}$, and $s_\emptyset = \emptyset$. Thus, for the $\{1,\ldots,D\}$-term of $P$ one has two possibilities: either this term is identically $\emptyset$ (if $P_{\{1,\ldots,D\}} = \emptyset$) or it equals $n^D$ (if $P_{\{1,\ldots,D\}} = \{\emptyset\}$).)

A set-polynomial $P$ has *degree* $d \in \mathbb{N}$, $\deg(P) = d$, if $P_\alpha = \emptyset$ when $|\alpha| > d$, $\alpha \subseteq \{1,\ldots,D\}$, and there exists $\alpha \subseteq \{1,\ldots,D\}$, $|\alpha| = d$, such that $P_\alpha \neq \emptyset$. *The empty set-polynomial*, that is, the set-polynomial all whose coefficients are empty, has, by definition, degree 0. Note that (since $D$ is fixed) all the set-polynomials have degrees not exceeding $D$.

*The term of degree $l$*, $l \leq D$, of a set-polynomial $P$ is the sum of its $\alpha$-terms for $|\alpha| = d$. *The constant term* of $P$ is its term of degree 0, that is, $P_\emptyset$, *the linear term* is the term of degree 1, *the leading term* of $P$ is its term of degree $\deg(P)$.

When possible, we will write set-polynomials in the more convenient form:
$$P(n) = \sum_{\alpha \subseteq \{1,\ldots,D\}} n^{|\alpha|} P_\alpha,$$

but this is not always possible as the example in (2.1) and examples 2, 3 below show.

*Examples.*
1. The general constant set-polynomial: $P(n) = P_\emptyset$, $P_\emptyset \in \mathcal{F}(V)$ (all the coefficients of $P$ with the exception of $P_\emptyset$ are empty).
2. The general linear set-polynomial:
$$P(n) = nP_{\{1\}} + \ldots + \{s : \ (s_1,\ldots,s_{i-1},s_{i+1},\ldots,s_D) \in P_{\{i\}}, \ s_i \in n\}$$
$$+ \ldots + P_{\{D\}}n + P_\emptyset,$$
where $P_{\{i\}} \in \mathcal{F}(S^{D-1})$.
3. The general quadratic set-polynomial:
$$P(n) = n^2 P_{\{1,2\}} + \ldots + P_{\{D-1,D\}}n^2 + nP_{\{1\}} + \ldots + P_{\{N\}}n + P_\emptyset.$$
(The dots contain the terms the writing of which is cumbersome.)
4. A set-polynomial of the maximal degree $D$:
$$P(n) = n^D + n^{D-1}P_{\{1,\ldots,D-1\}} + \ldots + nP_{\{1\}} + P_\emptyset.$$

2.4. *The sum $P + Q$ of two set-polynomials $P$, $Q$ is the set-polynomial, whose coefficients are the unions of the corresponding coefficients of $P$ and $Q$:*

$$(P+Q)_\alpha = P_\alpha + Q_\alpha.$$

The difference of set-polynomials is not always defined. We will say that a set-polynomial $Q$ is *dominated by* a set-polynomial $P$, $Q \leq P$, if for every $\alpha \subseteq \{1,\ldots,D\}$ one has $Q_\alpha \subseteq P_\alpha$. In this case we define the difference of $P$ and $Q$, $P - Q$, as the set-polynomial whose coefficients are the differences of the corresponding coefficients of $P$ and $Q$:

$$(P-Q)_\alpha = P_\alpha - Q_\alpha.$$

2.5. One more operation on set-polynomials is the "argument shifting": if $P$ is a set-polynomial and $m \in \mathcal{F}(S)$ then $R(n) = P(n+m)$, $n \in \mathcal{F}(S)$, is a set-polynomial as well; note, that $P \leq R$.

*Example.* $(n+m)^2 = n^2 + nm + mn + m^2$.

2.6. Given an element $s = (s_1,\ldots,s_d) \in S^d$, $d \in \mathbb{N}$, let $\pi_i(s) = s_i$. For a subset $a \in \mathcal{F}(S^d)$ let $\pi_i(a) = \bigcup_{s \in a} \pi_i(s)$. The support, supp$(a)$, of a set $a \in \mathcal{F}(S^d)$ is the subset of $S$ defined as

$$\mathrm{supp}(a) = \bigcup_{i=1}^d \pi_i(a).$$

The support, supp$(P)$, of a set-polynomial $P$ is the union of the supports of its coefficients:

$$\mathrm{supp}(P) = \bigcup_{\alpha \subseteq \{1,\ldots,D\}} \mathrm{supp}(P_\alpha).$$

2.7 LEMMA. *Let $P$ be a set-polynomial.*

1) *If $n, m \in \mathcal{F}(S)$ and $n \subseteq m$, then $P(n) \subseteq P(m)$.*

2) *If $n, m \in \mathcal{F}(S)$, then $P(n) \cup P(m) \subseteq P(n \cup m)$.*

3) *If $n \in \mathcal{F}(S)$, then $\mathrm{supp}\bigl(P(n)\bigr) \subseteq \mathrm{supp}(P) \cup n$.*

4) *If $m \in \mathcal{F}(S)$ and $R(n) = P(n+m)$, then $\mathrm{supp}(R) \subseteq \mathrm{supp}(P) \cup m$.*

5) *If $a \in \mathcal{F}(V)$ and $n \in \mathcal{F}(S)$, then $a \cap P(n) \subseteq P(\mathrm{supp}(a) \cap n)$.*

6) *If a set-polynomial $Q$ and a set $n \in \mathcal{F}(S)$ satisfy $n \cap \mathrm{supp}(Q) = \emptyset$, then for any $m \in \mathcal{F}(S)$ one has $P(n) \cap Q(m) \subseteq P(n \cap m)$.*

*Proof.* 1), 2), 3) and 4) are obvious.

5) Let $s = (s_1, \ldots, s_D) \in a \cap P(n)$. Let $\alpha \subseteq \{1, \ldots, D\}$ be such that $s$ belongs to the $\alpha$-term of $P$, that is, $s_{\bar{\alpha}} \in P_\alpha$ (where $\bar{\alpha} = \{1, \ldots, D\} \setminus \alpha$) and $s_i \in n$ for $i \in \alpha$. Then, since $s \in a$, one has $s_i \in \text{supp}(a)$, $i = 1, \ldots, D$. Hence, $s_i \in \text{supp}(a) \cap n$ for $i \in \alpha$, that is, $s \in P(\text{supp}(a) \cap n)$.

6) $P(n) \cap Q(m) \subseteq P(n \cap \text{supp}(Q(m))) \subseteq P(n \cap (\text{supp}(Q) \cup m))$
$= P(n \cap m)$. □

## 3. Actions of $\mathcal{F}(S)$

3.1. Let $W$ be a set.

Let us recall that *an action* of $\mathcal{F}(W)$ on a topological space $X$ is any mapping $T$ from $\mathcal{F}(W)$ into the set of continuous self-mappings of $X$, $a \mapsto T^a$, satisfying the following: for any $a, b \in \mathcal{F}(W)$, $a \cap b = \emptyset$, one has $T^{a \cup b} = T^a T^b$. It follows, in particular, that for any $a, b \in \mathcal{F}(W)$ one has $T^a T^b = T^b T^a$.

*Examples.*
1. Let $T$ be a continuous self-mapping of $X$. For $a \in \mathcal{F}(W)$, put $T^a = T^{|a|}$. Here, generally speaking, $T^{a \cup b} \neq T^a T^b$ when $a \cap b \neq \emptyset$.
2. Given a topological space $\mathbf{K}$, consider the product space $X = \mathbf{K}^{\mathcal{F}(W)} = \{\omega \colon \mathcal{F}(W) \longrightarrow \mathbf{K}\}$ with its natural product topology; $X$ is compact if $\mathbf{K}$ is. Define the action $T$ of $\mathcal{F}(W)$ on $X$ by

$$(T^a \omega)(b) = \omega(a \cup b) \text{ for } a, b \in \mathcal{F}(W), \, \omega \in X.$$

Contrary to the preceding example, $T^{a \cup b} = T^a T^b$ for any $a, b \in \mathcal{F}(W)$.

When $\mathbf{K} = \{1, \ldots, k\}$, $k \in \mathbb{N}$, we will call the topological space $\mathbf{K}^{\mathcal{F}(W)}$ *the space of $k$-colorings* of $\mathcal{F}(W)$ and denote it by $\Omega_{W,k}$, the corresponding action of $\mathcal{F}(W)$ on $\Omega_{W,k}$ will be denoted by $T_{W,k}$. If $W$ is countable then also $\mathcal{F}(W)$ is, and the topology on $\Omega_{W,k} = \{1, \ldots, k\}^{\mathcal{F}(W)}$ is given by some metric. In this metric, two colorings are close to one another if they coincide on all subsets of a large enough subset of $W$.

3.2. *A system* is a finite set of set-polynomials. *The support* of a system $A$ is the union of the supports of its elements: $\text{supp}(A) = \bigcup_{P \in A} \text{supp}(P)$.

Let $(X, \rho)$ be a metric space and $T$ be an action of $\mathcal{F}(V)$ on $X$. A system $A$ is *a system of recurrence with respect to $T$* if for any $\varepsilon > 0$ and any $H \in \mathcal{F}(S)$ there exists $N \in \mathcal{F}(S)$, $N \cap H = \emptyset$, such that for every $x \in X$ there exist $n \subseteq N$, $n \neq \emptyset$, and $a \subseteq \bigcup_{P \in A} P(N)$ for which $a \cap P(n) = \emptyset$ and $\rho(T^{a + P(n)} x, T^a x) < \varepsilon$ for any $P \in A$.

A system $A$ is *a system of recurrence* if $A$ is a system of recurrence with respect to any action of $\mathcal{F}(V)$ on any metric compact space.

A system $A$ is *a system of chromatic recurrence* if $A$ is a system of recurrence with respect to the action $T_{V,k}$ for any $k \in \mathbb{N}$. (It is clear that this definition does not depend on the concrete metric chosen on $\Omega_{V,k}$.)

3.3. The following proposition shows that the notion of a system of recurrence coincides with the notion of a system of chromatic recurrence. It gives also a clear combinatorial description of systems of recurrence.

PROPOSITION. *Let $A$ be a system. The following three statements are equivalent*:

(1) *$A$ is a system of recurrence.*
(2) *$A$ is a system of chromatic recurrence.*
(3) *For any $k \in \mathbb{N}$ and any $H \in \mathcal{F}(S)$ there exists $N \in \mathcal{F}(S)$, $N \cap H = \emptyset$, such that for any mapping $\omega\colon \mathcal{F}(V) \longrightarrow \{1,\ldots,k\}$ there exist $n \subseteq N$, $n \neq \emptyset$, and $a \subseteq \bigcup_{P \in A} P(N)$ such that, for each $P \in A$, one has $a \cap P(n) = \emptyset$ and $\omega\bigl(a + P(n)\bigr) = \omega(a)$.*

*Proof.*

$(1) \Longrightarrow (2)$ follows immediately from the definition.

$(2) \Longrightarrow (3)$ Let $k \in \mathbb{N}$ and $H \in \mathcal{F}(S)$. Since $A$ is a system of chromatic recurrence, there exists $N \in \mathcal{F}(S)$, $N \cap H = \emptyset$, such that for any point $\omega$ of the metric compact space $\Omega_{V,k}$ of $k$-colorings, $\omega\colon \mathcal{F}(V) \longrightarrow \{1,\ldots,k\}$, one can find sets $n \in \mathcal{F}(S)$, $n \neq \emptyset$, and $a \subseteq \bigcup_{P \in A} P(N)$ such that, for each $P \in A$, $a \cap P(n) = \emptyset$ and the distance between the points $T_{V,k}^{a+P(n)}\omega$ and $T_{V,k}^a \omega$ is as small as one needs to guarantee that their values at the empty set, $T_{V,k}^{a+P(n)}\omega(\emptyset) = \omega\bigl(a + P(n)\bigr)$ and $T_{V,k}^a\omega(\emptyset) = \omega(a)$, respectively, coincide.

$(3) \Longrightarrow (1)$ Let $(X, \rho)$ be a compact metric space, let $T$ be an action of $\mathcal{F}(V)$ on $X$, let $H \in \mathcal{F}(S)$ and let $\varepsilon > 0$. Fix an $\varepsilon/2$-net $\{x_1,\ldots,x_k\}$ in $X$.

By (3), there exist $N \in \mathcal{F}(S)$, $N \cap H = \emptyset$, such that for any mapping $\omega\colon \mathcal{F}(V) \longrightarrow \{1,\ldots,k\}$ there exist $n \subseteq N$, $n \neq \emptyset$, and $a \subseteq \bigcup_{P \in A} P(N)$ such that, for every $P \in A$, $a \cap P(n) = \emptyset$ and $\omega(a + n) = \omega(a)$.

Let $x \in X$; define $\omega_x\colon \mathcal{F}(V) \longrightarrow \{1,\ldots,k\}$ by the rule: for $b \in \mathcal{F}(V)$,

$$\omega_x(b) \text{ is any } j \text{ for which } \rho(T^b x, x_j) < \varepsilon/2.$$

Find $n \subseteq N$, $n \neq \emptyset$, and $a \subseteq \bigcup_{P \in A} P(N)$ such that $a \cap P(n) = \emptyset$ and $\omega_x(a + P(n)) = \omega_x(a)$ for any $P \in A$. The last equality means that for $j = \omega_x(a)$ one has

$$\rho(T^{a+P(n)} x, x_j) < \frac{\varepsilon}{2} \text{ and } \rho(T^a x, x_j) < \frac{\varepsilon}{2};$$

that is, $\rho(T^a x, T^{a+P(n)} x) < \varepsilon$. $\square$

3.4. Now we can state our main result, Theorem PHJ, in the following form:

THEOREM. *Let A be a system consisting of set-polynomials with empty constant terms (that is, $P(\emptyset) = \emptyset$ for every $P \in A$). Then A is a system of recurrence.*

The proof of this theorem will be given in Section 6.

3.5. In light of Proposition 3.3, we can also formulate Theorem 3.4 in the following equivalent "chromatic" form:

THEOREM. *Let A be a system consisting of set-polynomials whose constant terms are empty, let $H \in \mathcal{F}(S)$ and let $k \in \mathbb{N}$. Then there exists $N \in \mathcal{F}(S)$, $N \cap H = \emptyset$, such that for any mapping $\omega: \mathcal{F}(V) \longrightarrow \{1, \ldots, k\}$ there exist $n \subseteq N$, $n \neq \emptyset$, and $a \subseteq \bigcup_{P \in A} P(N)$ such that, for each $P \in A$, one has $a \cap P(n) = \emptyset$ and $\omega(a + P(n)) = \omega(a)$.*

*Remark.* Since for any $H \in \mathcal{F}(S)$ we can find $n \in \mathcal{F}(S)$ satisfying the statement of Theorem 3.5 and such that $n \cap H = \emptyset$, we have plenty of such $n$.

3.6. *Remark.* Note that in the proof of Proposition 3.3 we used only the fact that $X$ was totally bounded but neither the completeness of $X$ nor the continuity of the action $T$. This shows that Theorem 3.4 holds true as well if one considers totally bounded spaces and arbitrary (not necessarily continuous) actions on them. Theorem 3.5 being a corollary of Theorem 3.4 corresponding to actions on the space of finite colorings of an infinite set gives as a corollary a generalization of Theorem 3.4. One could ask, why, in such a case, we formulate this theorem for continuous actions on compact spaces instead of formulating it in the more general setup of arbitrary actions on totally bounded spaces. The reason is twofold. First, the noncontinuous version of Theorem 3.4 is a simple corollary of this theorem. Second, we follow a longstanding tradition.

3.7. Let us show now that Theorem 3.5 implies Theorem PHJ. To do this, it suffices to embed the parallelepiped $\{1, \ldots, N\}^d \times \{1, \ldots, q\}$, introduced in the formulation of Theorem PHJ, to the space $V = S^D$, where $D = d + 1$.

More exactly, let $r, d, q \in \mathbb{N}$ be given. Put $S = \mathbb{N}$, $D = d+1$, $V = S^D$, and define set-polynomials $P_0, P_1, \ldots, P_q$ by $P_0 = \emptyset$, $P_i(n) = n^d \times \{i\}$, $i = 1, \ldots, q$. Applying Theorem 3.5 to the system $A = \{P_0, \ldots, P_q\}$, find $N \in \mathbb{N}$ such that for any $r$-coloring of $\mathcal{F}(V)$ there exist a nonempty set $\gamma \subseteq \{1, \ldots, N\}$ and $a \subset \{1, \ldots, N\}^d \times \{1, \ldots, q\}$ such that $a \cap P_i(\gamma) = \emptyset$ for any $1 \leq i \leq q$ (and so, $a \cap (\gamma^d \times \{1 \ldots, q\}) = \emptyset$) and the set $\{a \cup P_i(\gamma) : i = 0, \ldots, q\}$ is monochromatic. Such $N$ satisfies the conclusion of Theorem PHJ.

The opposite implication, Theorem PHJ $\Longrightarrow$ Theorem 3.5 is also true. We will not prove it.

## 4. PET-induction

All the polynomials we will deal with from here on will have empty constant terms, and we will not mention this specifically. Throughout the sequel $S$ is a fixed infinite set, $D$ is a fixed positive integer and $V = S^D$.

4.1. We need some more definitions.

Two set-polynomials $P$ and $P'$ are called *equivalent* if their degrees coincide and their leading terms coincide too:

$P \sim P'$ if there is $d \leq D$ such that $P_\alpha = P'_\alpha = \emptyset$ for all $\alpha$ with $|\alpha| > d$,

and $P_\alpha = P'_\alpha$ for all $\alpha$ with $|\alpha| = d$.

Equivalent set-polynomials form equivalence classes; *the degree* of an equivalence class is the degree of any of its elements.

4.2. Attribute to each system $A$ its *weight vector* $w(A) = (w_1, \ldots, w_D)$, where the nonnegative integer $w_d$, $d = 1, \ldots, D$, denotes the number of equivalence classes of degree $d$ having a representative in $A$.

Define an order on the set of weight vectors as follows: $(w'_1, \ldots, w'_D) < (w_1, \ldots, w_D)$ if there exists $1 \leq d \leq D$ such that $w'_j = w_j$ for any $j > d$ and $w'_d < w_d$. We will say that *a system $A'$ precedes a system $A$ if $w(A') < w(A)$*.

4.3. We will prove Theorem 3.4 by induction on the well ordered set of weight vectors: we will be deducing that any system $A$ (consisting of set-polynomials with empty constant terms) is a system of recurrence from the assumption that all the systems preceding $A$ (and consisting of set-polynomials with empty constant terms) are systems of recurrence. This is the so-called *PET-induction* (PET stands for Polynomial Exhaustion Technique).

The beginning of the induction process is obvious: if the weight vector of a system $A$ is minimal, that is, equals $(0, \ldots, 0)$, $A$ may contain only constant and, so, empty set-polynomials. For such a system, the statement of Theorem 3.4 is trivial.

4.4. The main tools used in the course of the induction proof of PHJ are the following lemma and its corollary:

LEMMA.

1) Let $P$ be a set-polynomial, let $m \in \mathcal{F}(S)$ and let $R$ be the set-polynomial defined by $R(n) = P(n + m)$. Then $R \sim P$.
2) Let $P$, $P'$ and $Q \neq \emptyset$ be set-polynomials and let $P \sim P'$, $Q \leq P$ and $Q \leq P'$.
   a) If $P \not\sim Q$, then $P - Q \sim P' - Q$ and $\deg(P - Q) = \deg(P)$.

b) *If $P \sim Q$, then $\deg(P - Q) < \deg(P)$.*

We omit the proof.

4.5. COROLLARY.    *Let $A$ be a system.*

1) *If $A'$ consists only of set-polynomials of the form $R(n) = P(n+m)$, $P \in A$, $m \in \mathcal{F}(S)$, then $w(A') \leq w(A)$.*
2) *If $Q \in A$, $Q \neq \emptyset$, is such that $Q \leq P$ for every $P \in A$, then the system $A'' = \{P - Q, \ P \in A\}$ precedes $A$: $w(A'') < w(A)$.*

## 5. Systems having a minimal element

A set-polynomial is *a minimal element* of a system if it is dominated by all other elements of the system. Not every system has a minimal element: for example, the system $\{an, na\}$, where $a \in \mathcal{F}(S^{D-1})$, has no minimal element. The following three lemmas reduce the proof of Theorem 3.4 to the case of systems having a (nonempty) minimal element.

It will be convenient to assume in this section that all the set-polynomials under consideration do not reach the maximal degree $D$; we can always achieve this by increasing $D$ and embedding $S^D \longrightarrow S^{D+1}$, $s \mapsto s \times \{r\}$ for some $r \in S$.

5.1.    Let $A$ be a system. For every $d = 1, \ldots, D-1$, let $R_{d,1}, \ldots, R_{d,r_d}$ be all possible different terms of degree $d$, including the empty one as well, which occur in the set-polynomials of $A$. Choose, for each $d = 1, \ldots, D$, arbitrary elements $p_{d,i} \in S^{D-d}$, $i = 1, \ldots, r_d$, subject only to the condition that $\mathrm{supp}(\{p_{d,i}\})$ for all $d = 1, \ldots, D$, $i = 1, \ldots, r_d$, are pairwise disjoint. For each $P \in A$, $P = R_{1,i_1} + \ldots + R_{d,i_d}$, define a new set-polynomial $P'$:

$$P'(n) = n\{p_{1,i_1}\} + \ldots + n^d\{p_{d,i_d}\}.$$

Let $A'$ be the system consisting of these new set-polynomials: $A' = \{P' : P \in A\}$.

It is clear that the degrees of set-polynomials and their equivalence are preserved when we pass from $A$ to $A'$; hence, the weight vector of the system is preserved as well: $w(A') = w(A)$.

5.2. LEMMA.    *Under the notation of (5.1), $A$ is a system of recurrence if $A'$ is.*

*Proof.* Define a mapping $\psi \colon \mathcal{F}(V) \longrightarrow \mathcal{F}(V)$ by the following rules:

1) $\psi$ is additive, that is, $\psi(a + b) = \psi(a) + \psi(b)$ for any $a, b \in \mathcal{F}(V)$;

2) if a point $s' \in V$ has the form $(s_1, \ldots, s_d, p_{d,i})$ for some $1 \leq d \leq D-1$, $1 \leq i \leq r_d$, then

$$\psi(\{s'\}) = \bigcup_{\substack{\alpha \subseteq \{1,\ldots,D\} \\ |\alpha|=d}} \{s : s_{\bar{\alpha}} \in (R_{d,i})_\alpha, \ s_\alpha = (s_1, \ldots, s_d)\}$$

under the notation of 2.3;

3) $\psi(\{s'\}) = \{s'\}$ if $s'$ is not of the form $(s_1, \ldots, s_d, p_{d,i})$ for any $1 \leq d \leq D-1$, $1 \leq i \leq r_d$. (Since the elements $p_{d,i}$ have been chosen in such a way that their supports are pairwise disjoint, $\psi$ is well defined.)

It is easy to see that, for any $P \in A$ and any $n \in \mathcal{F}(S)$, one has

$$\psi(P'(n)) = P(n) \qquad (5.1)$$

and, for any $s' \in V$ and $s \in \psi(s')$, one has

$$\operatorname{supp}(\{s'\}) \setminus \operatorname{supp}(A') \subseteq \operatorname{supp}(\{s\}) \setminus \operatorname{supp}(A). \qquad (5.2)$$

Given a coloring $\omega \in \Omega_{V,k}$ one can define a new coloring $\psi(\omega)$ by $\psi(\omega) = \omega \circ \psi \in \Omega_{V,k}$. Let $k \in \mathbb{N}$ and $H \in \mathcal{F}(S)$ be given. Since $A'$ is a system of recurrence, by Proposition 3.3 there exists $N \in \mathcal{F}(S)$, $N \cap (H \cup \operatorname{supp}(A)) = \emptyset$, such that for every $\omega' \in \Omega_{V,k}$ there exist $n \subseteq N$, $n \neq \emptyset$, and $a' \subseteq \bigcup_{P' \in A'} P'(N)$ such that for any $P' \in A'$ one has $a' \cap P'(n) = \emptyset$ and $\omega'(a' + P'(n)) = \omega'(a')$.

Let $\omega \in \Omega_{V,k}$; find $n$ and $a'$ corresponding to $\omega' = \psi(\omega)$. Let $a = \psi(a')$. Then, by (5.1),

$$a \subseteq \psi\Big(\bigcup_{P' \in A'} P'(N)\Big) = \bigcup_{P \in A} P(N)$$

and

$$\omega(a + P(n)) = \omega(\psi(a' + P'(n))) = \omega'(a' + P'(n)) = \omega'(a') = \omega(a)$$

for any $P \in A$. It remains to prove that $a \cap P(n) = \emptyset$ for any $P \in A$. If this is the case, Proposition 3.3 shows that $A$ is a system of recurrence.

Let $s' = (s_1, \ldots, s_D) \in a'$, namely, $s_1, \ldots, s_d \in N$ and $(s_{d+1}, \ldots, s_D) \in P'_{\{1,\ldots,d\}}$ for some $1 \leq d \leq D-1$ and $P' \in A'$ under the notation of 2.3; let $s \in \psi(s')$. Since $a \cap P'(n) = \emptyset$, for some $1 \leq i \leq d$ one has $s_i \notin n$; that is,

$$\operatorname{supp}(\{s'\}) \cap (N \setminus n) \neq \emptyset.$$

Since $N \cap \operatorname{supp}(A) = \emptyset$, it follows from (5.2) that

$$\operatorname{supp}(\{s\}) \cap (N \setminus n) \neq \emptyset;$$

that is,
$$\operatorname{supp}(\{s\}) \not\subseteq n \cup \operatorname{supp}(A).$$

But, as $\operatorname{supp}\bigl(\bigcup_{P \in A} P(n)\bigr) \subseteq n \cup \operatorname{supp}(A)$, this means that $s \notin \bigcup_{P \in A} P(n)$. Since this is valid for any $s \in \psi(a')$, we obtain $a \cap \bigl(\bigcup_{P \in A} P(n)\bigr) = \emptyset$. □

5.3. Starting from $A'$, we will now construct still another system, $A''$, which will have a minimal element. We will show then that if $A''$ is a system of recurrence then $A'$ and hence (in light of Lemma 5.2) $A$ also are systems of recurrence.

Let $Q \in A'$ be an element of $A'$ of the minimal degree in $A'$. Put $A'' = \{P' + Q : P' \in A'\}$. It is clear that $w(A'') = w(A')$ and, furthermore, $Q$ is a minimal element of $A''$.

5.4. LEMMA.    *Under the notation of* 5.1 *and* 5.3, $A'$ *is a system of recurrence if* $A''$ *is.*

*Proof.* Let $Q(n) = n\{q_1\} + \ldots + n^e\{q_e\}$, $q_i \in S^{D-i}$, $i = 1, \ldots, e$, $e \leq D-1$. Define a mapping $\psi': \mathcal{F}(V) \longrightarrow \mathcal{F}(V)$ in the following way: the set $\psi'(b)$ is obtained from a set $b \in \mathcal{F}(V)$ by deleting some points. Namely, the point $s = (s_1, \ldots, s_D) \in b$ is deleted if and only if there exist a number $1 \leq d \leq e$, a set-polynomial $P = n\{p_1\} + \ldots + n^l\{p_l\} \in A'$ and a point $\tilde{s} \in b$, $\tilde{s} \neq s$, such that $(s_{d+1}, \ldots, s_D) = q_d$. That is $s = (s_1, \ldots, s_d, q_d)$, and $\tilde{s} = (s_1, \ldots, s_d, p_d)$. It is easy to see that, for any $P' \in A'$ and any $n \in \mathcal{F}(S)$, one has $\psi'\bigl(P'(n) + Q(n)\bigr) = P'(n)$.

Let $k \in \mathbb{N}$ and $H \in \mathcal{F}(S)$ be defined. Since $A''$ is a system of recurrence, by Proposition 3.3 there exists $N \in \mathcal{F}(S)$, $N \cap H = \emptyset$, such that for any $\omega' \in \Omega_{V,k}$ there exist $n \subseteq N$, $n \neq \emptyset$, and $a' \subseteq \bigcup_{P' \in A'} P'(N)$ such that for any $P' \in A'$ one has $a' \cap P'(n) = \emptyset$ and $\omega'\bigl(a' + P'(n) + Q(n)\bigr) = \omega'(a')$.

Let $\omega \in \Omega_{V,k}$; find $n$ and $a'$ corresponding to $\omega' = \psi'(\omega)$. Put $a = \psi'(a')$. Since $a \subseteq a'$, we have $a \subseteq \bigcup_{P' \in A'} P'(N)$ and $a \cap P'(n) = \emptyset$ for any $P' \in A'$. Furthermore, by the definitions of $\omega'$, $n$ and $a'$,
$$\omega(a) = \omega'(a') = \omega'\bigl(a' + P'(n) + Q(n)\bigr) = \omega\bigl(\psi'(a' + P'(n) + Q(n))\bigr)$$
for every $P' \in A'$.

Let us look at the set $\psi'\bigl(a' + P'(n) + Q(n)\bigr)$ more carefully. If $P' = Q$ this set is $\psi'(a' + Q(n))$. Let $P = n\{p_1\} + \ldots + n^l\{p_l\} \in A'$, let $d \in \{1, \ldots, e\}$, and let $s, \tilde{s} \in a' + Q(n)$, $s \neq \tilde{s}$, be such that $s = (s_1, \ldots, s_d, q_d)$, $\tilde{s} = (s_1, \ldots, s_d, p_d)$. By the choice of $A'$, one has $\tilde{s} \notin Q(n)$ (since either $p_d = q_d$ or $\operatorname{supp}(p_d) \cap \operatorname{supp}(q_i) = \emptyset$ for any $1 \leq i \leq e$), thus $\tilde{s} \in a'$. If it were that $s \in Q(n)$, that is $s_1, \ldots, s_d \in n$, one would have $\tilde{s} \in P(n)$ and so, $a' \cap P(n) \neq \emptyset$. Hence, $s, \tilde{s} \in a'$ and
$$\psi'(a' + Q(n)) = a + Q(n).$$

In the case $P' \neq Q$ a similar argument shows that

$$\psi'\big(a' + P'(n) + Q(n)\big) = \psi'(a') + \psi'\big(P'(n) + Q(n)\big) = a + P'(n).$$

This shows that $\omega(a) = \omega(a + P'(n))$ for any $P' \in A'$. Now, Proposition 3.3 applied once more gives that $A'$ is a system of recurrence. $\square$

## 6. Proof of Theorem 3.4

Let $A$ be a system consisting of set-polynomials with no constant terms, our goal is to prove that $A$ is a system of recurrence. Without loss of generality, we may and will assume that $A$ does not contain empty set-polynomials. Replace $A$ by the new system $A''$ built from $A$ as described in Section 5. This system has the same weight vector as $A$ and contains a nonempty minimal element. Lemma 5.2 and Lemma 5.4 show that it is enough to prove that $A''$ is a system of recurrence. We will proceed by PET-induction, described in Section 4; we will assume that the statement of Theorem 3.4 holds true for all systems preceding $A$ and, so, $A''$. Thus, it is enough to prove the following proposition.

PROPOSITION. *Let $A$ be a system consisting of set-polynomials with empty constant terms and containing a nonempty minimal element $Q \in A$: $Q \leq P$ for every $P \in A$. Furthermore, let all the systems preceding $A$ and consisting of set-polynomials without constant terms be systems of recurrence. Then $A$ is a system of recurrence as well.*

*Proof.* Let $H \in \mathcal{F}(S)$, let $(X, \rho)$ be a metric compact space and let $\varepsilon > 0$. Let $k \in \mathbb{N}$ be such that there are two points at a distance less than $\varepsilon/2$ among any $k+1$ points of $X$.

6.1. Put $H_0 = H \cup \mathrm{supp}(A)$, $\varepsilon_0 = \varepsilon/2k$ and define a new system $A_0$:

$$A_0 = \big\{R(n) = P(n) - Q(n) : P \in A\big\}.$$

By Corollary 4.5, we know that $A_0$ precedes $A$; so, $A_0$ is a system of recurrence and

$$\exists N_0 \in \mathcal{F}(S),\ N_0 \cap H_0 = \emptyset,\ \text{such that}$$

$$\forall x \in X\ \exists n \subseteq N_0\ \text{and}\ \exists a \subseteq \bigcup_{R \in A_0} R(N_0)\ \text{such that}\ n \neq \emptyset\ \text{and}$$

$$\forall R \in A_0\ \text{one has}\ a \cap R(n) = \emptyset\ \text{and}\ \rho(T^{a+R(n)}x, T^a x) < \varepsilon_0.$$

We put now $H_1 = H_0 \cup N_0$. Let $\varepsilon_1 > 0$ be such that the inequality $\rho(y_1, y_2) < \varepsilon_1$ implies $\rho(T^b y_1, T^b y_2) < \varepsilon/2k$ for any $b \subseteq \bigcup_{P \in A} P(N_0)$. Consider the system

$$A_1 = \big\{R(n) = P(n+m) - Q(n) - P(m) : P \in A,\ m \subseteq N_0\big\}.$$

By Corollary 4.5, $A_1$ precedes $A$; so,

$$\exists N_1 \in \mathcal{F}(S), \ N_1 \cap H_1 = \emptyset, \text{ such that}$$

$$\forall x \in X \ \exists n \subseteq N_1 \text{ and } \exists a \subseteq \bigcup_{R \in A_1} R(N_1) \text{ such that } n \neq \emptyset \text{ and}$$

$$\forall R \in A_1 \text{ one has } a \cap R(n) = \emptyset \text{ and } \rho(T^{a+R(n)}x, T^a x) < \varepsilon_1.$$

Continue this process: suppose that numbers $\varepsilon_0, \ldots, \varepsilon_i$ and sets $N_0, \ldots, N_i$ have been already chosen, and assume that the sets $H_0, \ldots, H_i$ and the systems $A_0, \ldots, A_i$ have been already defined. We put $H_{i+1} = H_i \cup N_i$. Let $\varepsilon_{i+1} > 0$ be such that the inequality $\rho(y_1, y_2) < \varepsilon_{i+1}$ implies $\rho(T^b y_1, T^b y_2) < \varepsilon/2k$ for any $b \subseteq \bigcup_{P \in A} P(N_0 \cup \ldots \cup N_i)$. Consider the system

$$A_{i+1} = \{R(n) = P(n+m) - Q(n) - P(m) : P \in A, \ m \subseteq N_0 \cup \ldots \cup N_i\}.$$

By Corollary 4.5, $A_{i+1}$ precedes $A$; so,

$$\exists N_{i+1} \in \mathcal{F}(S), \ N_{i+1} \cap H_{i+1} = \emptyset, \text{ such that}$$

$$\forall x \in X \ \exists n \subseteq N_{i+1}, \ a \subseteq \bigcup_{R \in A_{i+1}} R(N_{i+1}) \text{ such that } n \neq \emptyset \text{ and}$$

$$\forall R \in A_{i+1} \text{ one has } a \cap R(n) = \emptyset \text{ and } \rho(T^{a+R(n)}x, T^a x) < \varepsilon_{i+1}.$$

Continue the process of choosing $\varepsilon_i$, $N_i$ and defining $H_i$, $A_i$ up to $i = k$, and put $N = N_0 \cup \ldots \cup N_k$; we have $N \cap H = \emptyset$. □

6.2. Now fix an arbitrary point $x \in X$.
Find, applying the definition of $N_k$ to the point

$$y_k = T^{Q(N_k)}x,$$

sets $n_k \subseteq N_k$, $n_k \neq \emptyset$, and $a_k \subseteq \bigcup_{R \in A_k} R(N_k)$ such that for every $R \in A_k$ one has

$$a_k \cap R(n_k) = \emptyset \text{ and } \rho(T^{R(n_k)+a_k}y_k, T^{a_k}y_k) < \varepsilon_k.$$

Then, applying the definition of $N_{k-1}$ to the point

$$y_{k-1} = T^{Q(N_k)+Q(N_{k-1})+a_k}x,$$

find sets $n_{k-1} \subseteq N_{k-1}$, $n_{k-1} \neq \emptyset$, and $a_{k-1} \subseteq \bigcup_{R \in A_{k-1}} R(N_{k-1})$ such that for every $R \in A_{k-1}$ we have

$$a_{k-1} \cap R(n_{k-1}) = \emptyset \text{ and } \rho(T^{R(n_{k-1})+a_{k-1}}y_{k-1}, T^{a_{k-1}}y_{k-1}) < \varepsilon_{k-1}.$$

Continue this process: suppose that we have already found $n_k, \ldots, n_i$, $a_k, \ldots, a_i$. Applying the definition of $N_{i-1}$ to the point

$$y_{i-1} = T^{Q(N_k)+\ldots+Q(N_{i-1})+a_k+\ldots+a_i}x,$$

find sets $n_{i-1} \subseteq N_{i-1}$, $n_{i-1} \neq \emptyset$, and $a_{i-1} \subseteq \bigcup_{R \in A_{i-1}} R(N_{i-1})$ such that for every $R \in A_{i-1}$ one has

$$a_{i-1} \cap R(n_{i-1}) = \emptyset \text{ and } \rho(T^{R(n_{i-1})+a_{i-1}} y_{i-1}, T^{a_{i-1}} y_{i-1}) < \varepsilon_{i-1}.$$

Continue the process of choosing $n_i$, $a_i$ up to $i = 0$.

6.3. Note the following:

(i) Since $N_i \subseteq H_j$, $\operatorname{supp}(Q) \subseteq H_j$, $0 \leq i < j \leq k$, and $N_j \cap H_j = \emptyset$, by Lemma 2.7, one has $Q(N_i) \cap Q(N_j) = Q(\emptyset) = \emptyset$ for $i \neq j$, $i, j = 0, \ldots, k$.

(ii) By definition,

$$a_i \subseteq \bigcup_{R \in A_i} R(N_i) \subseteq \bigcup_{\substack{P \in A \\ m \subseteq N_0 \cup \ldots \cup N_{i-1}}} \bigl(P(N_i + m) - P(m) - Q(N_i)\bigr), \quad 0 \leq i \leq k.$$

This gives $a_i \cap Q(N_i) = \emptyset$, $i = 0, \ldots, k$.

(iii) Furthermore, $\operatorname{supp}(a_i) \subseteq H_{i+1} = H \cup N_0 \cup \ldots \cup N_i$, $0 \leq i \leq k$. Since $N_j \cap H_{i+1} = \emptyset$, $0 \leq i < j \leq k$, we have $a_i \cap Q(N_j) = \emptyset$, $0 \leq i < j \leq k$.
Since $a_j \subseteq \bigcup_{R \in A_j} R(N_j)$ and $N_j \cap (N_i \cup \operatorname{supp}(Q)) = \emptyset$, $0 \leq i < j \leq k$, we have also $a_j \cap Q(N_i) = \emptyset$, $0 \leq i < j \leq k$.

(iv) Let $P \in A$, let $0 \leq i < j \leq k$ and $0 \leq l \leq k$. Then, by Lemma 2.7,

$$P(n_j + \ldots + n_{i+1}) \cap Q(N_l) \subseteq \begin{cases} \emptyset, & l \leq i \text{ or } l > j \\ Q(n_l), & i < l \leq j \end{cases}$$

and, hence,

$$P(n_j + \ldots + n_{i+1}) \cap \bigl(Q(N_l) - Q(n_l)\bigr) = \emptyset.$$

(v) Under the notation of (iv), using Lemma 2.7 and taking into account that $\operatorname{supp}(a_l) \subseteq H \cup N_0 \cup \ldots \cup N_l$, one has

$$P(n_j + \ldots + n_{i+1}) \cap a_l = P\bigl((n_j + \ldots + n_{i+1}) \cap \operatorname{supp}(a_l)\bigr) \cap a_l$$
$$= \begin{cases} \emptyset, & l \leq i \\ P(n_l + \ldots + n_{i+1}) \cap a_l, & l > i. \end{cases}$$

In the case $l > i$, represent $P(n_l + \ldots + n_{i+1})$ in the form

$$\bigl(P(n_l + \ldots + n_{i+1}) - P(n_{l-1} + \ldots + n_{i+1}) - Q(n_l)\bigr) + P(n_{l-1} + \ldots + n_{i+1}) + Q(n_l).$$

Since $P(n_l + \ldots + n_{i+1}) - P(n_{l-1} + \ldots + n_{i+1}) - Q(n_l) \in A_l$, we have

$$\bigl(P(n_l + \ldots + n_{i+1}) - P(n_{l-1} + \ldots + n_{i+1}) - Q(n_l)\bigr) \cap a_l = \emptyset.$$

Furthermore, since $Q(n_l) \subseteq Q(N_l)$, by (iii) one has $Q(n_l) \cap a_l = \emptyset$, and we obtain

$$P(n_j + \ldots + n_{i+1}) \cap a_l \subseteq P(n_{l-1} + \ldots + n_{i+1}) \cap a_l.$$

Recall that
$$a_l \subseteq \bigcup_{\substack{P' \in A \\ m \subseteq N_0 \cup \ldots \cup N_{l-1}}} \left(P'(N_l + m) - P'(m) - Q(N_l)\right).$$

Now, by Lemma 2.7, for any $P' \in A$ and any $m \subseteq N_0 \cup \ldots \cup N_{l-1}$, one has
$$P(n_{l-1} + \ldots + n_{i+1}) \cap P'(N_l + m) \subseteq P'(m),$$
and, hence,
$$P(n_j + \ldots + n_{i+1}) \cap a_l = \emptyset.$$

6.4. Define points $x_i$, $i = 0, \ldots, k$, by
$$x_i = T^{a_k + \ldots + a_0 + Q(N_k) + \ldots + Q(N_0) - Q(n_i) - \ldots - Q(n_0)} x.$$

We are going to show that, for any $0 \le i \le j \le k$ and for any $P \in A$,
$$(6.1) \qquad \rho(T^{P(n_j + \ldots + n_{i+1})} x_j, x_i) \le \frac{\varepsilon}{2k}(j - i).$$

We will prove this by induction on $j - i$. When $j = i$ the statement is trivial. (We remind the agreement from subsection 1.1.4: the indices of the entries of expressions like "$n_j + \ldots + n_{i+1}$ for $i \le j$" are assumed to be decreasing; when $i = j$ this expression is assumed to be vacuous and equal to zero.) We will derive the validity of (6.1) for $i, j$, where $i < j$, from its validity for $i, j - 1$.

By the definition of $n_j$,
$$\rho(T^{a_j + P(n_j + \ldots + n_{i+1}) - P(n_{j-1} + \ldots + n_{i+1}) - Q(n_j)} y_j, T^{a_j} y_j) < \varepsilon_j,$$

where
$$y_j = T^{Q(N_k) + \ldots + Q(N_j) + a_k + \ldots + a_{j+1}} x.$$

So, by the choice of $\varepsilon_j$,
$$\rho(T^{P(n_{j-1} + \ldots + n_{i+1}) + a_{j-1} + \ldots + a_0 + Q(N_{j-1}) + \ldots + Q(N_0) - Q(n_{j-1}) - \ldots - Q(n_0)}$$
$$T^{a_j + P(n_j + \ldots + n_{i+1}) - P(n_{j-1} + \ldots + n_{i+1}) - Q(n_j)} y_j,$$
$$T^{P(n_{j-1} + \ldots + n_{i+1}) + a_{j-1} + \ldots + a_0 + Q(N_{j-1}) + \ldots + Q(N_0) - Q(n_{j-1}) - \ldots - Q(n_0)} T^{a_j} y_j) < \frac{\varepsilon}{2k}.$$

Using the additivity of $T$, the definition of $y_j$, $x_j$ and the facts established in 6.3, we see that
$$T^{P(n_{j-1} + \ldots + n_{i+1}) + a_{j-1} + \ldots + a_0 + Q(N_{j-1}) + \ldots + Q(N_0) - Q(n_{j-1}) - \ldots - Q(n_0)}$$
$$T^{a_j + P(n_j + \ldots + n_{i+1}) - P(n_{j-1} + \ldots + n_{i+1}) - Q(n_j)} y_j$$
$$= T^{P(n_j + \ldots + n_{i+1}) + a_k + \ldots + a_0 + Q(N_k) + \ldots + Q(N_0) - Q(n_j) - \ldots - Q(n_0)} x$$
$$= T^{P(n_j + \ldots + n_{i+1})} x_j$$

and
$$T^{P(n_{j-1}+\ldots+n_{i+1})+a_{j-1}+\ldots+a_0+Q(N_{j-1})+\ldots+Q(N_0)-Q(n_{j-1})-\ldots-Q(n_0)}T^{a_j}y_j$$
$$= T^{P(n_{j-1}+\ldots+n_{i+1})+a_k+\ldots+a_0+Q(N_k)+\ldots+Q(N_0)-Q(n_{j-1})-\ldots-Q(n_0)}x$$
$$= T^{P(n_{j-1}+\ldots+n_{i+1})}x_{j-1}.$$

Since, by the induction hypothesis,
$$\rho(T^{P(n_{j-1}+\ldots+n_{i+1})}x_{j-1}, x_i) \le \frac{\varepsilon}{2k}(j-i-1),$$
we obtain (6.1).

6.5. By the choice of $k$, among the $k+1$ points $x_0, \ldots, x_k$ there are two, say $x_i, x_j$, $0 \le i < j \le k$, for which $\rho(x_i, x_j) < \varepsilon/2$. Put
$$n = n_j + \ldots + n_{i+1}, \ a = a_k + \ldots + a_0 + Q(N_k) + \ldots + Q(N_0) - Q(n_j) - \ldots - Q(n_0).$$
Then $x_j = T^a x$ and, for every $P \in A$, from 6.3 one has $a \cap P(n) = \emptyset$; hence,
$$\rho(T^{a+P(n)}x, T^a x) = \rho(T^{P(n)}x_j, x_j) \le \rho(T^{P(n)}x_j, x_i) + \rho(x_j, x_i)$$
$$< \frac{\varepsilon}{2k}(j-i) + \frac{\varepsilon}{2} \le \varepsilon.$$
Furthermore, $n \subseteq N$, $n \ne \emptyset$ and $a \subseteq \bigcup_{P \in A} P(N)$. □

## 7. Applications

In this section we will derive from Theorem PHJ some general combinatorial results, as well as corollaries pertaining to topological recurrence. As a matter of fact, the combinatorial corollaries that we deal with (even if they look as remote from dynamics as the van der Waerden Theorem or the geometric Ramsey Theorem (see subsections 0.2 and 0.3)) are also instances of recurrence, namely of *chromatic recurrence*.

This section is organized in the following way. We bring first the definitions of (families of) topological and chromatic recurrence. We produce next a corollary of PHJ, Proposition 7.3, which will serve as a bridge between the set-theoretical recurrence inherent in PHJ and chromatic and topological recurrence related to commutative semigroups and actions thereof. The rest of the section is devoted to derivation of special cases of interest (which include the results mentioned in the introduction).

7.1. Let $G$ be a commutative semigroup. We say that an indexed family $\mathcal{R} = \{R_w, \ w \in \mathcal{W}\}$ of nonempty subsets of $G$, $R_w \subseteq G$ for $w \in \mathcal{W}$, is *a family of (topological) recurrence* if for any action $T = \{T^g, \ g \in G\}$ of $G$ on a metric compact space $(X, \rho)$ by continuous self-mappings and any $\varepsilon > 0$ there exist $x \in X$ and $w \in \mathcal{W}$ such that $\rho(x, T^g x) < \varepsilon$ for every $g \in R_w$. In the special case when all the members of the family of recurrence $\mathcal{R}$ are one-element sets, $R_w = \{g_w\}$, we will say that $\{g_w, w \in \mathcal{W}\}$ is *a set of recurrence*. Similarly, a family $\mathcal{R} = \{R_w, w \in \mathcal{W}\}$ is *a family of chromatic recurrence* if for any finite coloring of $G$ there is $w \in \mathcal{W}$ and an element $h \in G$ such that the set $h + R_w$ is monochromatic. One can show that a family $R$ is a family of recurrence if and only if $R$ is a family of chromatic recurrence (see 3.3). Note also that if $G$ is a subgroup of a commutative group $G'$, then $\mathcal{R}$ is a family of recurrence in $G$ if and only if $\mathcal{R}$ is a family of recurrence in $G'$. Some examples of families of recurrence were already brought in subsections 0.13–0.17.

7.2. The following simple fact will be used in the sequel.

LEMMA. *Given a homomorphism $\varphi: G_1 \longrightarrow G_2$ of commutative semigroups and a family of recurrence $\{R_w\}_{w \in \mathcal{W}}$ in $G_1$, the family $\{\varphi(R_w)\}_{w \in \mathcal{W}}$ is a family of recurrence in $G_2$.*

*Proof.* Indeed, given a finite coloring $\chi: G_2 \longrightarrow \{1, \ldots, r\}$ of $G_2$, it induces a coloring $\chi'$ of $G_1$ by $\chi' = \chi \circ \varphi$, and if $w \in W$ and $h_1 \in G_1$ are such that the set $h_1 + R_w$ is monochromatic with respect to $\chi'$, then its image $\varphi(h_1) + \varphi(R_w)$ in $G_2$ is monochromatic with respect to $\chi$. $\square$

7.3. The following proposition is modeled on Corollary 3.2 in [BM1]. Despite its cumbersome appearance it is just the right result enabling one to obtain general combinatorial corollaries from PHJ. An additional reason for appreciation of this proposition is that its special cases or modifications turn out to serve as a very convenient (and so far, the only) tool in obtaining strong density results generalizing the classical Szemerédi Theorem on arithmetic progressions in various polynomial directions (see [BM1], [BM2], [L]). In what follows, we denote by $G[z_1, z_2, \ldots]$ the additive semigroup of formal polynomials in noncommuting variables $z_1, z_2, \ldots$ with coefficients from a commutative semigroup $G$.

PROPOSITION. *Let $r, q, n \in \mathbb{N}$, let $G$ be a commutative semigroup and let $p_1, \ldots, p_q \in G[x_1, \ldots, x_n]$ be polynomials without constant terms. There exists $N \in \mathbb{N}$ and a finite set $\mathcal{Q} \subset G[y_{1,1}, \ldots, y_{1,N}, y_{2,1}, \ldots, y_{2,N}, \ldots, y_{n,1}, \ldots, y_{n,N}]$ of polynomials without constant terms such that for any $r$-coloring of $\mathcal{Q}$ there exist a nonempty $\gamma \subseteq \{1, \ldots, N\}$ and a polynomial $h \in \mathcal{Q}$ for which the polynomials*

$$h(y_{1,1}, \ldots, y_{n,N}) + p_t\Big(\sum_{j \in \gamma} y_{1,j}, \ldots, \sum_{j \in \gamma} y_{n,j}\Big), \quad t = 1, \ldots, q,$$

*are all in $\mathcal{Q}$ and have the same color.*

*Proof.* Let $d = \max\{\deg p_1, \ldots, \deg p_q\}$ and let $N = N(r, d, q)$ be the bound guaranteed by Theorem PHJ. Put $V = \{1, \ldots, N\}^d$.

For every point $v = (j_1, \ldots, j_d)$ of $V$ and every monomial $\eta = cx_1^{d_1} \ldots x_n^{d_n}$, where $c \in G$ and $d_1, \ldots, d_n$ are nonnegative integers satisfying $1 \leq d_1 + \ldots + d_n \leq d$, define a polynomial $\Phi_{v,\eta}$ of $nN$ variables $y_{1,1}, \ldots, y_{n,N}$ by

$$\Phi_{v,\eta}(y_{1,1}, \ldots, y_{n,N}) = c\Big(\prod_{k=1}^{d_1} y_{1,j_k}\Big)\Big(\prod_{k=d_1+1}^{d_1+d_2} y_{2,j_k}\Big) \cdots \Big(\prod_{k=d_1+\ldots+d_{n-1}+1}^{d_1+\ldots+d_n} y_{n,j_k}\Big)$$

if $j_{d_1+\ldots+d_n} = j_{d_1+\ldots+d_n+1} = \ldots = j_d$, and $\Phi_{v,\eta}(y_{1,1}, \ldots, y_{n,N}) = 0$ otherwise. Expand $\Phi$ to $\mathcal{F}(V)$ by $\Phi_{a,\eta} = \sum_{v \in a} \Phi_{v,\eta}$ for $a \in \mathcal{F}(V)$ and a monomial $\eta$. It is easy to see that for any nonempty subset $\gamma \subseteq \{1, \ldots, N\}$ we have

$$\Phi_{\gamma^d,\eta}(y_{1,1}, \ldots, y_{n,N}) = \eta\Big(\sum_{j \in \gamma} y_{1,j}, \ldots, \sum_{j \in \gamma} y_{q,j}\Big).$$

Now, for any polynomial $p(x_1, \ldots, x_n) = \sum_l \eta_l$ define $\Phi_{a,p} = \sum_l \Phi_{a,\eta_l}$. We obtain a mapping $\Phi: \mathcal{F}(V) \times G[x_1, \ldots, x_n] \longrightarrow G[y_{1,1}, \ldots, y_{n,N}]$, additive with respect to the first argument in the following sense: for any polynomial $p$ and for any $a, b \in \mathcal{F}(V)$ with $a \cap b = \emptyset$ we have $\Phi_{a \cup b,p} = \Phi_{a,p} + \Phi_{b,p}$. In addition, for any nonempty subset $\gamma \subseteq \{1, \ldots, N\}$,

$$\Phi_{\gamma^d,p}(y_{1,1}, \ldots, y_{n,N}) = p\Big(\sum_{j \in \gamma} y_{1,j}, \ldots, \sum_{j \in \gamma} y_{q,j}\Big).$$

Define a mapping $\varphi: \mathcal{F}(V)^q \longrightarrow G[y_{1,1}, \ldots, y_{n,N}]$ by $\varphi(a_1, \ldots, a_q) = \sum_{t=1}^{q} \Phi_{a_t,p_t}$, and take as $\mathcal{Q}$ the image of $\varphi$. Now, given an $r$-coloring $\chi$ of $\mathcal{Q}$, it defines the coloring $\chi' = \chi \circ \varphi$ of $\mathcal{F}(V)^q$. Find $a \in \mathcal{F}(V)^q$ and a nonempty subset $\gamma \subseteq \{1, \ldots, N\}$ such that $a \cap (\gamma^d \times \{1, \ldots, q\}) = \emptyset$ and the sets $a$, $a \cup (\gamma^d \times \{1\})$, $\ldots$, $a \cup (\gamma^d \times \{q\})$ have the same color, and put $h = \varphi(a)$. Since $\varphi(a \cup (\gamma^d \times \{t\})) = h(y_{1,1}, \ldots, y_{n,N}) + p_t(\sum_{j \in \gamma} y_{1,j}, \ldots, \sum_{j \in \gamma} y_{n,j})$, $t = 1, \ldots, q$, such $\gamma$ and $h$ satisfy the conclusion of the proposition. □

7.4. Under the definition in subsection 7.1, Proposition 7.3 provides us with the following corollary.

COROLLARY. *Let $G$ be a commutative semigroup and let $p_1, \ldots, p_q$ be polynomials of $n$ variables over $G$ without constant terms. Then*

$$\left\{\left\{p_1\Big(\sum_{j \in \gamma} y_{1,j}, \ldots, \sum_{j \in \gamma} y_{n,j}\Big), \ldots, p_q\Big(\sum_{j \in \gamma} y_{1,j}, \ldots, \sum_{j \in \gamma} y_{n,j}\Big)\right\} : \gamma \subset \mathbb{N},\ 0 < |\gamma| < \infty\right\}$$

*is a family of recurrence (in the semigroup of polynomials $G[y_{m,j},\ m = 1, \ldots, n,\ j = 1, 2, \ldots]$).*

7.5. Let now $G$ be a (not necessarily commutative) ring and let $p_1, \ldots, p_q$ be polynomials of $n$ variables over $G$ without constant terms. Let $g_{m,j} \in G$, $m = 1, \ldots, n$, $j \in \mathbb{N}$. Then we can define a homomorphism of $G$-algebras $\varphi\colon G[y_{m,j},\ m=1,\ldots,n,\ j=1,2,\ldots] \longrightarrow G$ by putting $y_{m,j} \stackrel{\varphi}{\mapsto} g_{m,j}$. Applying Corollary 7.4 and Lemma 7.2, we obtain:

PROPOSITION. *The family*

$$\left\{ \left\{ p_1\Big(\sum_{j\in\gamma} g_{1,j},\ldots,\sum_{j\in\gamma} g_{n,j}\Big),\ldots,p_q\Big(\sum_{j\in\gamma} g_{1,j},\ldots,\sum_{j\in\gamma} g_{n,j}\Big) \right\} : \gamma \subset \mathbb{N},\ 0 < |\gamma| < \infty \right\}$$

*is a family of recurrence in $G$.*

7.6. As an example, taking the polynomials $p_1 = y_1 z_1, \ldots, p_q = y_q z_q$ we obtain the following statement:

*Given sequences $\{g_{1,j}\}_{j\in\mathbb{N}}, \{h_{1,j}\}_{j\in\mathbb{N}}, \ldots, \{g_{q,j}\}_{j\in\mathbb{N}}, \{h_{q,j}\}_{j\in\mathbb{N}}$ of elements of $G$, the family*

$$\left\{ \left\{ \Big(\sum_{j\in\gamma} g_{1,j}\Big)\Big(\sum_{j\in\gamma} h_{1,j}\Big),\ldots,\Big(\sum_{j\in\gamma} g_{q,j}\Big)\Big(\sum_{j\in\gamma} h_{q,j}\Big) \right\} : \gamma \subset \mathbb{N},\ 0 < |\gamma| < \infty \right\}$$

*is a family of recurrence in $G$.*

To see that this statement is indeed Proposition 0.15 in disguise, one should recall that a family $\mathcal{R}$ is a family of recurrence if and only if it is a family of chromatic recurrence (see subsection 7.1).

7.7. As a matter of fact, we even did not need the ring structure on $G$ in 7.5. Let $G$ be the tensor product, $G = G_1 \otimes_\mathbb{N} \ldots \otimes_\mathbb{N} G_D$, of some commutative semigroups $G_i$, $i = 1, \ldots, D$ (that is, the commutative semigroup generated by $g_1 \otimes \ldots \otimes g_D$, $g_i \in G_i$, with relations $g_1 \otimes \ldots \otimes (g_i' + g_i'') \otimes \ldots \otimes g_D = g_1 \otimes \ldots \otimes g_i' \otimes \ldots \otimes g_D + g_1 \otimes \ldots \otimes g_i'' \otimes \ldots \otimes g_D$, $i = 1, \ldots, D$). Then the conclusion of Proposition 7.5 holds true if we replace products of elements of $G$ by tensor products of elements of corresponding $G_i$. (It follows from the fact that $\mathbb{N}[y_1, \ldots, y_D] \simeq \mathbb{N}[y_1] \otimes_\mathbb{N} \ldots \otimes_\mathbb{N} \mathbb{N}[y_D]$ as additive semigroups.) Since the general statement is too cumbersome, we will confine ourselves to an instructive example. Namely, the example of subsection 7.6 takes under this generalization the following form:

Let $G = G_1 \otimes_\mathbb{N} G_2$, where $G_1$, $G_2$ are commutative semigroups. Given sequences $\{g_{1,j}\}_{j\in\mathbb{N}}, \ldots, \{g_{q,j}\}_{j\in\mathbb{N}}$ in $G_1$ and $\{h_{1,j}\}_{j\in\mathbb{N}}, \ldots, \{h_{q,j}\}_{j\in\mathbb{N}}$ in $G_2$, the family

$$\left\{ \left\{ \Big(\sum_{j\in\gamma} g_{1,j}\Big) \otimes \Big(\sum_{j\in\gamma} h_{1,j}\Big),\ldots,\Big(\sum_{j\in\gamma} g_{q,j}\Big) \otimes \Big(\sum_{j\in\gamma} h_{q,j}\Big) \right\} : \gamma \subset \mathbb{N},\ 0 < |\gamma| < \infty \right\}$$

*is a family of recurrence in $G$.*

7.8.   Let us utilize this statement in the case $G_1 = (\mathbb{N}, \cdot)$, $G_2 = (\mathbb{N}, +)$. Then $G = G_1 \otimes_{\mathbb{N}} G_2 \simeq G_1$, where the isomorphism is given by $u \otimes v \mapsto u^v$. We obtain the following version of Proposition 0.16:

Let $\{u_{1,j}\}_{j \in \mathbb{N}}, \{v_{1,j}\}_{j \in \mathbb{N}}, \ldots, \{u_{q,j}\}_{j \in \mathbb{N}}, \{v_{q,j}\}_{j \in \mathbb{N}}$ be sequences of positive integers. Then

$$\left\{ \left\{ (\prod_{j \in \gamma} u_{1,j})^{\sum_{j \in \gamma} v_{1,j}}, \ldots, (\prod_{j \in \gamma} u_{q,j})^{\sum_{j \in \gamma} v_{q,j}} \right\} : \gamma \subset \mathbb{N},\ 0 < |\gamma| < \infty \right\}$$

is a family of recurrence in $(\mathbb{N}, \cdot)$.

Again, to see that this statement is the same as Proposition 0.16, one uses the fact that a family $\mathcal{R}$ is a family of recurrence if and only if it is a family of chromatic recurrence.

7.9.   Let us now formulate a multidimensional corollary of Proposition 7.3.

PROPOSITION.     Let $n, t \in \mathbb{N}$. Let $G$ be a ring, let $P: G^n \longrightarrow G^t$ be a polynomial mapping over $G$ with $P(0) = 0$, let $\{g_j\}_{j \in \mathbb{N}}$ be a sequence of elements of $G^n$ and let $F$ be a finite configuration in $G^n$. For any finite coloring of $G^t$ there are $h \in G^t$ and a nonempty finite set $\gamma \subset \mathbb{N}$ such that for $(x_1, \ldots, x_n) = \sum_{j \in \gamma} g_j \in G^n$, the set

$$\{h + P(x_1 v_1, \ldots, x_n v_n) : (v_1, \ldots, v_n) \in F\}$$

is monochromatic in $G^t$.

7.10.   To prove Proposition 7.9, we first derive from Corollary 7.4 a multidimensional extension of Proposition 7.5. Again, let $G$ be a ring and let $\{g_{m,j}\}_{j \in \mathbb{N}}$, $m = 1, \ldots, n$, be $n$ sequences of elements of $G$. Let $t \in \mathbb{N}$; putting $\varphi(y_{l,m,j}) = g_{m,j}$, $m = 1, \ldots, n$, $j \in \mathbb{N}$, $l = 1, \ldots, t$, we obtain a homomorphism of $G$-modules

$$\varphi: \bigoplus_{l=1}^{t} G[y_{l,m,j},\ m = 1, \ldots, n,\ j = 1, 2, \ldots] \longrightarrow G^t.$$

Let $p_{i,l}$, $i = 1, \ldots, q$, $l = 1, \ldots, t$, be polynomials of $n$ variables over $G$ taking on zero at zero. Since $\bigoplus_{l=1}^{t} G[y_{l,m,j},\ m = 1, \ldots, n,\ j = 1, 2, \ldots]$ is a subgroup of the group $G[y_{l,m,j},\ m = 1,\ldots,n,\ j = 1, 2,\ldots,\ l = 1,\ldots,t]$, the family

$$R = \left\{ \left\{ \left(p_{1,1}\left(\sum_{j\in\gamma} y_{1,1,j},\ldots,\sum_{j\in\gamma} y_{1,n,j}\right),\ldots,p_{1,t}\left(\sum_{j\in\gamma} y_{t,1,j},\ldots,\sum_{j\in\gamma} y_{t,n,j}\right)\right),\ldots, \right.\right.$$

$$\left.\left.\left(p_{q,1}\left(\sum_{j\in\gamma} y_{1,1,j},\ldots,\sum_{j\in\gamma} y_{1,n,j}\right),\ldots,p_{q,t}\left(\sum_{j\in\gamma} y_{t,1,j},\ldots,\sum_{j\in\gamma} y_{t,n,j}\right)\right)\right\} : \gamma \subset \mathbb{N},\ 0 < |\gamma| < \infty \right\}$$

is a family of recurrence in $\bigoplus_{l=1}^{t} G[y_{l,m,j},\ m = 1,\ldots,n,\ j = 1, 2,\ldots]$. By Lemma 7.2, its image $\varphi(R) \subseteq G^t$ is also a family of recurrence. That is, we have shown the following fact:

LEMMA. *The family*

$$\left\{ \left\{ \left(p_{1,1}\left(\sum_{j\in\gamma} g_{1,j},\ldots,\sum_{j\in\gamma} g_{n,j}\right),\ldots,p_{1,t}\left(\sum_{j\in\gamma} g_{1,j},\ldots,\sum_{j\in\gamma} g_{n,j}\right)\right),\ldots,\right.\right.$$

$$\left.\left.\left(p_{q,1}\left(\sum_{j\in\gamma} g_{1,j},\ldots,\sum_{j\in\gamma} g_{n,j}\right),\ldots,p_{q,t}\left(\sum_{j\in\gamma} g_{1,j},\ldots,\sum_{j\in\gamma} g_{n,j}\right)\right)\right\} : \gamma \subset \mathbb{N},\ 0 < |\gamma| < \infty \right\}$$

*is a family of recurrence in $G^t$.*

Now, to prove Proposition 7.9 it is enough to apply this lemma to polynomials $p_{i,v}(x_1,\ldots,x_n) = \pi_i\bigl(P(x_1 v_1,\ldots,x_n v_n)\bigr)$, $i = 1,\ldots,t$, $v = (v_1,\ldots,v_n) \in F$, where $\pi_i$ is the projection of $G^t$ onto its $i$-th component.

7.11. In fact, the statement of Proposition 7.9 is trivial if the ring $G$ is finite. Indeed, in this case we always can find a set $\gamma \subset \mathbb{N}$ so that $\sum_{j\in\gamma} g_j$ vanishes in $G^n$. (See, however, Proposition 7.13 below.) Let now $G$ be an infinite ring; then we can choose $g_j \in G$, $j \in \mathbb{N}$, so that $\sum_{j\in\gamma} g_j \neq 0$ for every nonempty finite set $\gamma \subset \mathbb{N}$. Then Proposition 7.9 gives us immediately the following corollary (which has Theorem 0.8 as a special case).

COROLLARY. *Let $n, t \in \mathbb{N}$, let $G$ be an infinite ring, let $P: G^n \longrightarrow G^t$ be a polynomial mapping over $G$ with $P(0) = 0$, and let $F$ be a finite configuration in $G^n$. For any finite coloring of $G^t$ there are $h \in G^t$ and a nonzero $x \in G$ such that the set $\{h + P(xF)\}$ is monochromatic.*

7.12. Moreover, if $G$ is an infinite commutative integral domain and $P$ separates elements of $F$ (that is, $|P(F)| = |F|$), then we can choose the sequence $g_1, g_2, \ldots \in G$ such that for every $x \in G$ of the form $x = \sum_{j\in\gamma} g_j$, where

$\gamma$ is a nonempty finite subset of $\mathbb{N}$, one has $P(xv_1) \neq P(xv_2)$ for all distinct $v_1, v_2 \in F$. Indeed, the polynomial equation $P\big((g+x)v_1\big) = P\big((g+x)v_2\big)$ has only finitely many solutions (if any) with respect to $g$ if $x \in G$ and $v_1, v_2 \in F$ are fixed. Thus, if $g_1, \ldots, g_k$ have been already chosen, we can pick $g_{k+1} \in G$ such that $P\big((g_{k+1}+x)v_1\big) \neq P\big((g_{k+1}+x)v_2\big)$ for all pairs $v_1 \neq v_2 \in F$ and all $x$ of the form $x = \sum_{j \in \gamma} g_j$, $\gamma \subseteq \{1, \ldots, k\}$.

COROLLARY. *Let $n, t \in \mathbb{N}$, let $G$ be an infinite commutative integral domain, let $P \colon G^n \longrightarrow G^t$ be a polynomial mapping with $P(0) = 0$, and let $F$ be a finite configuration in $G^n$ such that $P$ is one-to-one on $F$. For any finite coloring of $G^t$ there are $h \in G^t$ and $x \in G$ such that the set $\{h + P(xF)\}$ is of the same cardinality as $F$ and monochromatic.*

7.13. We bring for completeness a "finitary" version of Corollary 7.12. We omit the proof, which may be given by utilizing Proposition 7.3.

PROPOSITION. *For any $n, t, d, a, r \in \mathbb{N}$ there is $N \in \mathbb{N}$ such that, given an integral domain $G$ of cardinality $\geq N$, a set $F \subset G^n$ of cardinality $\leq a$, a polynomial mapping $P \colon G^n \longrightarrow G^t$ of degree $\leq d$, satisfying $P(0) = 0$ and which is one-to-one on $F$, and an $r$-coloring of $G^t$, one can find $h \in G^t$ and $x \in G$ for which the set $\{h + P(xF)\}$ is monochromatic and has the same cardinality as $F$.*

7.14. Recall that, given a sequence of elements $\{u_j\}_{j \in \mathbb{N}}$ of a commutative semigroup $G$, *the* IP-*set* generated by $\{u_j\}$ is defined as the set of all finite sums of elements of $\{u_j\}$ with distinct indices,

$$FS(\{u_j\}) = \Big\{\sum_{j \in \gamma} u_j : \gamma \subset \mathbb{N},\ 0 < |\gamma| < \infty\Big\}.$$

Given a sequence of vectors $\{(g_{1,j}, \ldots, g_{n,j})\}_{j \in \mathbb{N}}$ in $G^n$, in Lemma 7.10 put $g_m = \sum_{j \in \gamma} g_{m,j}$, $m = 1, \ldots, n$; then

$$(g_1, \ldots, g_n) \in FS\big(\{(g_{1,j}, \ldots, g_{n,j})\}\big).$$

The sets having non-empty intersection with every IP-set are called IP$^*$-*sets*; we have, consequently:

COROLLARY. *Let $P \colon G^n \longrightarrow G^t$ be a polynomial mapping with $P(0) = 0$. Given a finite set $F \subset G^n$ and a finite coloring $\chi$ of $G^t$, $\chi \colon G^t \longrightarrow \{1, \ldots, r\}$, one can find a color $l$, $1 \leq l \leq r$, such that the vectors $(x_1, \ldots, x_n)$ in $G^n$ for which there is $h \in G^t$ with the property that $\chi\big(h + P(x_1v_1, \ldots, x_nv_n)\big) = l$ for all $(v_1, \ldots, v_n) \in F$, form an IP$^*$-set.*

7.15. We conclude this section by formulating a general topological multiple (and multiparameter) recurrence result which follows directly from Lemma 7.10. It generalizes Theorem C in [BL]. Let $T_1, \ldots, T_t$ be a set of pairwise commuting actions of the additive group of a ring $G$ by (continuous) self-mappings on a (pre-)compact metric space $(X, \rho)$. Then $G^t$ acts on $X$ by $T^{(g_1,\ldots,g_t)} = T_1^{g_1} \ldots T_t^{g_t}$. We get:

COROLLARY. *Let $p_{i,l}$, $i = 1, \ldots, q$, $l = 1, \ldots, t$, be polynomials of $n$ variables over $G$ taking on zero at zero. For any $\varepsilon > 0$ the set of vectors $(g_1, \ldots, g_n) \in G^n$ for which there exists $x \in X$ such that for every $i = 1, \ldots, q$ one has $\rho\left(\prod_{l=1}^{t} T_l^{p_{i,l}(g_1,\ldots,g_n)} x, x\right) < \varepsilon$, is an IP*-set.*

## 8. Appendix: An abstract version of PHJ

In the main part of the paper we introduced and treated the notion of a set-polynomial. We considered also some exotic "polynomial" expressions, like $\pi_\alpha^{\sigma_\alpha}$ (cf. 0.16). The natural question arises, in what sense are these objects polynomials? And what is the polynomial in general? This question would not look difficult had we dealt with rings: any polynomial is the sum of monomials. But structures which we have to deal with are not always provided with a multiplication, and consequently monomials are not defined there. As a matter of fact the proofs of our propositions are based on the following property of polynomials: every polynomial $P$ vanishes after an application of sufficiently many "differentiations" of the form $(D_m P)(n) = P(n+m) - P(n)$. We are now going to convert this property of polynomials into a definition.

8.1. Given a set $S$, the set of its finite subsets $\mathcal{F}(S)$ is a commutative semigroup with respect to the operation of taking the unions. (As a matter of fact, we are mainly concerned with unions of disjoint sets.)

We call a mapping $P$ from $\mathcal{F}(S)$ into a commutative semigroup $G$ a *polynomial mapping of degree* 0 if $P$ is constant, *a polynomial mapping of degree $\leq d \in \mathbb{N}$* if for any $m \in \mathcal{F}(S)$ there exists a polynomial mapping $D_m P : \mathcal{F}(S \setminus m) \longrightarrow G$ of degree $\leq d-1$ such that $P(n \cup m) = P(n) + (D_m P)(n)$ for every $n \in \mathcal{F}(S)$, $n \cap m = \emptyset$. We say that $P$ is *a polynomial mapping of degree $d$* if it is polynomial of degree $\leq d$ and not of degree $\leq d - 1$.

8.2. It is easy to give a natural representation for any polynomial mapping of degree $\leq 1$. Let $S$ be a set and let $P : \mathcal{F}(S) \longrightarrow G$ be a polynomial mapping with $\deg P \leq 1$. For every $s \in S$, let $\varphi(\{s\}) \in G$ be such that $P(n \cup \{s\}) = P(n) + \varphi(s)$ for all $n \in \mathcal{F}(S)$, $s \notin n$. Put $\varphi(\emptyset) = P(\emptyset)$. Then, by induction on $|n|$,

(8.1) $$P(n) = \varphi(\emptyset) + \sum_{s \in n} \varphi(\{s\}) \text{ for all } n \in \mathcal{F}(S).$$

On the other hand, given any mapping $\varphi \colon \{\{s\} : s \in S\} \cup \{\emptyset\} \longrightarrow G$, the mapping $P \colon \mathcal{F}(S) \longrightarrow G$ defined by formula (8.1) is polynomial of degree $\leq 1$.

8.3. We will now characterize polynomial mappings of higher degrees. Let $S$ be a set and let $d \in \mathbb{N}$; define $S^{\leq d} = \{a \in \mathcal{F}(S) : |a| \leq d\}$.

THEOREM. *Let $G$ be a commutative semigroup. For any mapping $\varphi \colon S^{\leq d} \longrightarrow G$ the mapping $P \colon \mathcal{F}(S) \longrightarrow G$ defined by*

(8.2) $$P(n) = \sum_{\substack{a \subseteq n \\ |a| \leq d}} \varphi(a)$$

*is polynomial of degree $\leq d$, and any polynomial mapping of degree $\leq d$ can be obtained in this fashion with the help of a suitable mapping $\varphi$.*

*Proof.* Let $\varphi \colon S^{\leq d} \longrightarrow G$ be given, let $P$ be defined by (8.2). Let $m \in \mathcal{F}(S)$; define $D_m P \colon \mathcal{F}(S \setminus m) \longrightarrow G$ by $(D_m P)(n) = \sum_{\substack{a \subseteq n \cup m \\ a \cap m \neq \emptyset \\ |a| \leq d}} \varphi(a)$, $n \in \mathcal{F}(S)$, $n \cap m = \emptyset$. Then

$$P(n \cup m) = \sum_{\substack{a \subseteq n \cup m \\ |a| \leq d}} \varphi(a) = \sum_{\substack{a \subseteq n \\ |a| \leq d}} \varphi(a) + \sum_{\substack{a \subseteq n \cup m \\ a \cap m \neq \emptyset \\ |a| \leq d}} \varphi(a) = P(n) + (D_m P)(n),$$

and $(D_m P)(n) = \sum_{\substack{b \subseteq n \\ |b| \leq d-1}} \psi(b)$, where we put $\psi(b) = \sum_{\substack{c \subseteq m \\ |c| \leq d - |b|}} \varphi(b \cup c)$ for $b \in \mathcal{F}(S \setminus m)$. By induction on $d$ we may conclude that the mapping $D_m P$ is polynomial of degree $\leq d - 1$ and hence $P$ is polynomial of degree $\leq d$.

In the opposite direction, let $P \colon \mathcal{F}(S) \longrightarrow G$ be a polynomial mapping of degree $\leq d$. Introduce an arbitrary linear ordering on $S$. Define $\varphi \colon S^{\leq d} \longrightarrow G$ in the following way. Put $\varphi(\emptyset) = P(\emptyset)$. Let $a \in \mathcal{F}(S)$, $1 \leq |a| \leq d$, and let $s$ be the minimal element of $a$. Let $D_s P \colon \mathcal{F}(S) \longrightarrow G$ be a polynomial mapping of degree $\leq d - 1$ such that $P(n \cup \{s\}) = P(n) + (D_s P)(n)$ for all $n \in \mathcal{F}(S \setminus \{s\})$. By induction on $d$, there exists a mapping $\varphi_s \colon (S \setminus \{s\})^{\leq d-1} \longrightarrow G$ such that $(D_s P)(n) = \sum_{\substack{b \subseteq n \\ |b| \subseteq d-1}} \varphi_s(b)$ for all $n \in \mathcal{F}(S)$, $s \notin n$. Put $\varphi(a) = \varphi_s(a \setminus \{s\})$.

Check now that (8.2) holds for every $n \in \mathcal{F}(S)$. It is so for $n = \emptyset$. Let $n \in \mathcal{F}(S)$, $n \neq \emptyset$, and $s$ be the minimal element of $n$. Then

$$P(n) = P(n \setminus \{s\}) + (D_s P)(n \setminus \{s\}) = P(n \setminus \{s\}) + \sum_{\substack{b \subseteq n \setminus \{s\} \\ |b| \leq d-1}} \varphi_s(b)$$

$$= P(n \setminus \{s\}) + \sum_{\substack{b \subseteq n \setminus \{s\} \\ |b| \leq d-1}} \varphi(b \cup \{s\}).$$

Since by induction on $|n|$ we may assume that $P(n \setminus \{s\}) = \sum_{\substack{a \subseteq n \setminus \{s\} \\ |a| \leq d}} \varphi(a)$, this implies (8.2). □

8.4. Given a set $S$ and a commutative semigroup $G$, define *a homomorphism* $\varphi \colon \mathcal{F}(S) \longrightarrow G$ as a mapping satisfying the equation $\varphi(n \cup m) = \varphi(n) + \varphi(m)$ for all $n, m \in \mathcal{F}(S)$ with $n \cap m = \emptyset$. Note that any nonconstant homomorphism is a polynomial mapping of degree 1 and is defined by its values at one-element subsets of $S$ by $\varphi(n) = \sum_{s \in n} \varphi(\{s\})$. Theorem 8.3 now says that any polynomial mapping $P \colon \mathcal{F}(S) \longrightarrow G$ of degree $\leq d$ can be represented as the composition $P = \varphi \circ \mathcal{P}^{(d)}$, where $\mathcal{P}^{(d)} \colon \mathcal{F}(S) \longrightarrow \mathcal{F}(S^{\leq d})$ is "the universal polynomial of degree $d$" defined by

$$\mathcal{P}^{(d)}(n) = \{a \subseteq n : |a| \leq d\},$$

and $\varphi \colon \mathcal{F}(S^{\leq d}) \longrightarrow G$ is a homomorphism.

8.5. Note now that, under the definitions of this section, in the formulation of PHJ we dealt with a polynomial mapping $\mathcal{P}^d \colon \mathcal{F}(S) \longrightarrow \mathcal{F}(S^d)$ (where $S^d$ was the set of $d$-tuples of $S$) defined by $\mathcal{P}^d(n) = n^d$, $n \in \mathcal{F}(S)$ (see Theorems PHJ and PHJt, where $S = \{1, \ldots, N\}$). Introduce an arbitrary linear ordering on $S$ and embed $S^{\leq d}_* = S^{\leq d} \setminus \{\emptyset\}$ into $S^d$ by

$$\{s_1, \ldots, s_k\} \mapsto (s_1, \ldots, s_k, s_k, \ldots, s_k)$$

where $s_1 < \ldots < s_k$. Then $\mathcal{P}^{(d)}(n) \setminus \{\emptyset\} = \mathcal{P}^d(n) \cap S^{\leq d}_*$ for all $n \in \mathcal{F}(S)$. It follows that every polynomial mapping $P \colon \mathcal{F}(S) \longrightarrow G$ of degree $\leq d$ with $P(\emptyset) = 0$ (under the assumption that $G$ has a zero element) can be obtained as the composition of $\mathcal{P}^d$ and a homomorphism $\varphi \colon \mathcal{F}(S^d) \longrightarrow G$. To see this, it suffices to take the homomorphism $\varphi \colon \mathcal{F}(S^{\leq d}) \longrightarrow G$ whose existence is guaranteed by Theorem 8.3, and extend $\varphi$ to $\mathcal{F}(S^d)$ by $\varphi(a) = \varphi(a \cap S^{\leq d})$. Thus, the polynomial mapping $\mathcal{P}^d$ can be considered as "the most general" polynomial mapping of degree $d$ without constant term. This shows that in PHJ we did not lose generality.

8.6. We can now bring a general group-theoretical corollary of Theorem PHJ in the following form:

THEOREM. *Let $S$ be an infinite set, let $F$ be a finite set and let $P$ be a polynomial mapping from $\mathcal{F}(S \times F)$ into a commutative group $G$ with $P(\emptyset) = 0$. Then for any finite coloring of $G$ there are $h \in G$ and a nonempty $\gamma \in \mathcal{F}(S)$ such that the set*
$$\{h + P(\gamma \times c) : c \subseteq F\}$$
*is monochromatic.*

*Proof.* Let $\deg P = d$. For every $c \subseteq F$ the mapping $P_c \colon \mathcal{F}(S) \longrightarrow G$, $P_c(n) = P(n \times c)$, is polynomial of degree $\leq d$ and hence is the composition of $\mathcal{P}^d$ and some homomorphism $\varphi_c \colon \mathcal{F}(S^d) \longrightarrow G$. Define a homomorphism $\varphi$ from $\mathcal{F}(S^d)^{\{c,\ c \subseteq F\}} = \mathcal{F}\bigl(S^d \times \{c,\ c \subseteq F\}\bigr)$ into $G$ as the sum of $\varphi_c$, $c \subseteq F$. A coloring of $G$ now induces a coloring of $\mathcal{F}\bigl(S^d \times \{c,\ c \subseteq F\}\bigr)$; let $a \in \mathcal{F}\bigl(S^d \times \{c,\ c \subseteq F\}\bigr)$ and $\gamma \in \mathcal{F}(S)$ be those given for this coloring by Theorem PHJ. Then $h = \varphi(a)$ and $\gamma$ satisfy the conclusion of the theorem. $\square$

8.7. Let $G$ and $G'$ be commutative groups. We say that $p \colon G \longrightarrow G'$ is *a polynomial mapping of degree $0$* if $p$ is constant, $p$ is *a polynomial mapping of degree $\leq d$*, $d \in \mathbb{N}$, if for any $a \in G$, the mapping $D_a p \colon G \longrightarrow G'$ defined by $(D_a p)(b) = p(a+b) - p(b)$ is polynomial of degree $\leq d-1$.

Note that if $G$ is a ring and $P \colon G^n \longrightarrow G^t$ is a conventional polynomial mapping (that is, $P = (p_1, \ldots, p_t)$ where $p_1, \ldots, p_t$ are polynomials over $G$ in $n$ variables), then $P$ as a mapping of the additive group of $G^n$ into the additive group of $G^t$ is polynomial.

8.8. We can now formulate a corollary of PHJ which generalizes the propositions of Section 7.

COROLLARY. *Let $G$ and $G'$ be commutative groups and let $p_1, \ldots, p_k \colon G \longrightarrow G'$ be polynomial mappings with $p_1(0) = \ldots = p_k(0) = 0$. For any finite coloring $\chi$ of $G'$, $\chi \colon G' \longrightarrow \{1, \ldots, r\}$, there is an $l$, $1 \leq l \leq r$, such that the set*
$$R = \bigl\{g \in G : \text{there is } h \in G' \text{ such that}$$
$$\chi\bigl(h + p_1(g)\bigr) = \ldots = \chi\bigl(h + p_k(g)\bigr) = l\bigr\}$$
*is an $\mathrm{IP}^*$-set in $G$.*

*Proof.* Let $\{g_l\}$ be any sequence in $G$. Define a mapping $P \colon \mathcal{F}\bigl(\mathbb{N} \times \{1, \ldots, k\}\bigr) \longrightarrow G$ by
$$P(a) = \sum_{i=1}^{k} p_i\Bigl(\sum_{l \in a_i} g_l\Bigr),$$
where $a = (a_1, \ldots, a_k) \in \mathcal{F}(\mathbb{N})^k = \mathcal{F}\bigl(\mathbb{N} \times \{1, \ldots, k\}\bigr)$. It can be checked directly from the definitions that $P$ is a polynomial mapping.

Let now $\chi$ be a finite coloring of $G$. Pick $h \in G$ and a nonempty set $\gamma \in \mathcal{F}(\mathbb{N})$ for which the set $S = \{h + P(\gamma \times \{i\}), \ i = 1, \ldots, k\}$ is monochromatic. Then for $g = \sum_{l \in \gamma} g_l$ we have $S = \{h + p_i(g), \ i = 1, \ldots, k\}$, that is, $g \in R$. It means that $R$ intersects $FS(\{g_l\})$ nontrivially, and thus is an IP$^*$-set. $\square$

The Ohio State University, Columbus, OH
*E-mail addresses*: vitaly@math.ohio-state.edu
leibman@math.ohio-state.edu